\documentclass[twocolumn,letterpaper]{IEEEAerospaceCLS}  % only supports two-column, letterpaper format

\usepackage[]{graphicx}    % We use this package in this document
\newcommand{\ignore}[1]{}  % {} empty inside = %% comment

\usepackage{graphbox}
\usepackage{array}
\usepackage{booktabs}
\usepackage{amssymb}
\usepackage[version=4]{mhchem}
\usepackage{siunitx}
\usepackage{longtable,tabularx}
\usepackage{booktabs}
\usepackage{mathtools}
\usepackage{multirow}
\usepackage{caption}
\usepackage{subcaption}

\usepackage{amsthm}

\usepackage[T1]{fontenc}\usepackage{tikz}
\usepackage{amsmath, amsfonts}
\usetikzlibrary{arrows.meta}

\usepackage[linesnumbered,ruled,vlined,norelsize,noend]{algorithm2e}
\SetInd{5pt}{10pt}
\SetNlSty{relax}{}{\enspace}

\newcommand{\bu}{\boldsymbol{u}}
\newcommand{\bx}{\boldsymbol{x}}

\usepackage{xcolor}
\definecolor{light-gray}{gray}{0.95}
\newcommand{\code}[1]{\colorbox{light-gray}{\texttt{#1}}}

\usepackage{hyperref}

\usepackage[font=small,skip=4pt]{caption}
\begin{document}
\title{Towards Robust Spacecraft Trajectory Optimization \\ via Transformers}

\author{%
Yuji Takubo$^*$\\ 
Dept. of Aeronautics \& Astronautics\\
Stanford University\\
496 Lomita Mall, Stanford, CA 94305\\
ytakubo@stanford.edu
\and
Tommaso Guffanti$^*$\\ 
Dept. of Aeronautics \& Astronautics\\
Stanford University\\
496 Lomita Mall, Stanford, CA 94305\\
tommaso@stanford.edu
\and 
Daniele Gammelli$^*$\\
Dept. of Aeronautics \& Astronautics\\
Stanford University\\
496 Lomita Mall, Stanford, CA 94305\\
gammelli@stanford.edu
\and 
Marco Pavone \\ 
Dept. of Aeronautics \& Astronautics\\
Stanford University\\
496 Lomita Mall, Stanford, CA 94305\\
pavone@stanford.edu
\and
Simone D'Amico \\ 
Dept. of Aeronautics \& Astronautics\\
Stanford University\\
496 Lomita Mall, Stanford, CA 94305\\
damicos@stanford.edu
%%%% IMPORTANT: Use the correct copyright information--IEEE, Crown, or U.S. government. %%%%%
\thanks{\footnotesize $^*$ The first three authors have contributed equally to the paper.}
% \thanks{$^\dagger$ The project’s website can be found at:  \url{https://robust-art.github.io/}}
%\def\thefootnote{*}\footnotetext{Equal contribution.}
% }
% \thanks{\footnotesize 979-8-3503-0462-6/24/$\$31.00$ \copyright2024 IEEE}              % This creates the copyright info that is the correct 2024 data.
%\thanks{{U.S. Government work not protected by U.S. copyright}}         % Use this copyright notice only if you are employed by the U.S. Government.
%\thanks{{979-8-3503-0462-6/24/$\$31.00$ \copyright2024 Crown}}          % Use this copyright notice only if you are employed by a crown government (e.g., Canada, UK, Australia).
%\thanks{{979-8-3503-0462-6/24/$\$31.00$ \copyright2024 European Union}}    % Use this copyright notice is you are employed by the European Union.
% }
}

\maketitle

\thispagestyle{plain}
\pagestyle{plain}

\maketitle

\thispagestyle{plain}
\pagestyle{plain}

\begin{abstract}
Future multi-spacecraft missions require robust autonomous trajectory optimization capabilities to ensure safe and efficient rendezvous operations. 
This capability hinges on solving non-convex optimal control problems in real-time, although traditional iterative methods such as sequential convex programming impose significant computational challenges. 
To mitigate this burden, the Autonomous Rendezvous Transformer (ART)~\cite{art_ieeeaero24} introduced a generative model trained to provide near-optimal initial guesses.
This approach provides convergence to better local optima (e.g., fuel optimality), improves feasibility rates, and results in faster convergence speed of optimization algorithms through warm-starting.
This work extends the capabilities of ART to address robust chance-constrained optimal control problems. 
Specifically, ART is applied to challenging rendezvous scenarios in Low Earth Orbit (LEO), ensuring fault-tolerant behavior under uncertainty. 
Through extensive experimentation, the proposed warm-starting strategy is shown to consistently produce high-quality reference trajectories, achieving up to 30\% cost improvement and 50\% reduction in infeasible cases compared to conventional methods, demonstrating robust performance across multiple state representations.
Additionally, a post hoc evaluation framework is proposed to assess the quality of generated trajectories and mitigate runtime failures, marking an initial step toward the reliable deployment of AI-driven solutions in safety-critical autonomous systems such as spacecraft.

\noindent Project website and Code: \url{https://robust-art.github.io/}
\end{abstract} 

\tableofcontents

%%%%%%%%%%%%%%%%%%%%%%%%%%%%%%%%%%%%%%
\section{Introduction}
%%%%%%%%%%%%%%%%%%%%%%%%%%%%%%%%%%%%%%
The new paradigm of in-space servicing, assembly, and logistics operations around the Earth and the Moon will require advanced autonomous rendezvous capabilities \cite{NASA_ARTEMIS} \cite{zea2024orbital}. 
Within this context, onboard trajectory optimization algorithms play a critical role, enabling spacecraft to compute sequences of states and controls that optimize mission objectives while adhering to operational constraints \cite{scp_2000} \cite{liu_scp_2014} \cite{malyuta_scp_2022}. 
Among other approaches, the integration of Machine Learning (ML) techniques as a “warm-start” for trajectory optimization algorithms offers a promising path to enhance traditional trajectory generation while ensuring strict constraint satisfaction, which is an essential requirement for the deployment of safety-critical systems \cite{art_ieeeaero24} \cite{art_ral24} \cite{Banerjee_2020} \cite{cauligi2021coco}. 
Notably, recent works have shown that Transformer networks \cite{chen2021decision} \cite{janner2021offline} \cite{art_ieeeaero24} \cite{art_ral24} \cite{art_acc25} are well-suited for this task, treating physical trajectories as tokenized sequences of states, controls, and, optionally, performance metrics. 
This approach has demonstrated that Transformer-generated warm-starts improve sequential optimizers by achieving more fuel-efficient trajectories, requiring fewer optimizer iterations, and reducing computation time, making them an appealing solution for spacecraft onboard operations. 

This paper builds upon the prior works presented in \cite{art_ieeeaero24} \cite{art_ral24} to address the remaining challenges in deploying Transformer-assisted trajectory optimization onboard spacecraft. 
This is approached across two key dimensions.
The first dimension focuses on enhancing robustness guarantees, specifically: (i) integrating a chance-constraint formulation \cite{boyd_chance} in the presence of uncertainty sources relevant within space operations (e.g., stemming from navigation, actuation, and unmodeled system dynamics), and (ii) integrating fault-tolerant constraint formulations \cite{guffanti_jgcd_2023} \cite{berning2024chance} \cite{marsillach2020fail} (e.g., ensuring passive safety in the event of control loss or actuation malfunctions). 
The second dimension evaluates the feasibility and benefits of deploying Transformer-based trajectory optimization algorithms onboard spacecraft, considering task performance, computational efficiency, and operational practicality. 
Through this, this paper offers four contributions to the current state-of-the-art:
\begin{enumerate}
    \item This work demonstrates how enforcing the aforementioned robustness guarantees requires tailored design choices in both the trajectory sequence definition and the warm-starting framework. Specifically, this work introduces novel performance parameter representations to encapsulate relevant information regarding robustness requirements. These tailored representations enable the Transformer model to generate trajectories that are robust to uncertainties and potential faults. 
    % allowing cognizant trajectory generation by the Transformer;
    \item This work empirically validates the proposed framework in challenging rendezvous scenarios in Low Earth Orbit, accounting for uncertainty and ensuring fault tolerance.
    \item This work demonstrates how Transformer models can generate high-performance rendezvous trajectories across different state representations, such as Cartesian \cite{hill1878researches} and Relative Orbit Elements \cite{damico_phd_2010}, for relative motion between spacecraft.
    \item This work addresses operational considerations for deploying Transformer-based trajectory optimization onboard spacecraft. Particularly, a novel post hoc evaluation scheme is introduced to ensure the quality of the model's outputs by detecting out-of-distribution outputs.
\end{enumerate}

\subsubsection{Notation} 
$\mathbb{R}$ and ${\mathbb{R}}_{\geq0}$ denote the space of real numbers and nonnegative real numbers, respectively; 
$\mathbb{N}$ denotes the set of natural numbers;
$\Pr(\cdot)$ denotes the probability operator; 
$\lVert\cdot\rVert_p$ denotes the $p$-norm; 
$\mathrm{Diag}\{\cdot\}$ denotes a diagonal matrix; 
For a positive semi-definite matrix $A\succeq0$, $A^{1/2}$ is a square-root of $A$ and satisfies $A = A^{1/2} \left(A^{1/2}\right)^\top$.

% %%%%%%%%%%%%%%%%%%%%%%%%%%%%%%%%%%%%%%%%%%%%%
% \section{Preliminaries}\label{sec_II}
% %%%%%%%%%%%%%%%%%%%%%%%%%%%%%%%%%%%%%%%%%%%%%

% This section provides a brief overview of the fundamental concepts of trajectory generation using trasformers for warm-starting sequential optimization relevant for this paper.

% \subsection{Transformers for Trajectory Generation}

% \subsection{Warm-starting Trajectory Optimization Problems}

%%%%%%%%%%%%%%%%%%%%%%%%%%%%%%%%%%%%%%%%%%%%%
\section{Robust Fault-tolerant Spacecraft Trajectory Optimization}\label{sec_II}
%%%%%%%%%%%%%%%%%%%%%%%%%%%%%%%%%%%%%%%%%%%%%

This work aims to solve an Optimal Control Problem (OCP) that computes an optimal maneuver sequence for a spacecraft, guiding it from an initial to a terminal state while accounting for uncertainties.
These uncertainties arise from various sources such as navigation errors, actuation errors, and unmodeled systems dynamics. 
In addition, ensuring safety during contingencies is crucial in spacecraft rendezvous. 
Specifically, a sufficient separation between spacecraft must be maintained not only along the nominal (controlled) trajectory but also in situations where control is unexpectedly lost at least for a certain duration---this concept is referred to as passive safety or drift safety.
To address this complex guidance problem, that accounts for both uncertainty and contingency, the following Chance-Constrained Optimal Control Problem (CC-OCP) \cite{blackmore2011chance} is formulated for a discrete-time linear dynamical system: 
\begin{equation} \label{OCP}
\resizebox{0.98\linewidth}{!}{$
\begin{aligned}
& \underset{\boldsymbol{x}(t_i), \boldsymbol{u}(t_i), \boldsymbol{X}(t_i))}{\text{minimize}} & & J = \sum_{i=1}^{N} {\mathcal{J}}(\boldsymbol{x}(t_i), \boldsymbol{u}(t_i), t_i) &\\
& \text{subject to} & & \boldsymbol{x}(t_{i+1}) = \boldsymbol{\Phi}(t_{i+1}, t_{i}) \Big(\boldsymbol{x}(t_{i}) + \boldsymbol{B}(t_{i}) \boldsymbol{u}(t_{i}) \Big)  & \forall i \in [1, N) \\
&&& \boldsymbol{X}(t_{i}) = f(\boldsymbol{x}(t_i), \boldsymbol{u}(t_i)) & \forall i \in [1, N] \\
&&& \Pr ({\mathcal{M}}^{*}(\boldsymbol{x}(t_i), \boldsymbol{X}(t_i)) \leq 0 ) \geq 1-\Delta_i & \forall i \in [1, N] \\
&&& \boldsymbol{u}(t_i) \in \mathcal{U} & \forall i \in [1, N] \\
&&&  \boldsymbol{x}(t_1) =  \boldsymbol{x}_{0} \hspace{0.2cm} \textrm{and} \hspace{0.2cm} \boldsymbol{x}(t_N) =  \boldsymbol{x}_{f}, & \\
\end{aligned}$}
\end{equation}
where the state is characterized by the mean $\boldsymbol{x} \in { {\mathbb{R}}^{6}}$ and error covariance $\boldsymbol{X} \in {{\mathbb{R}}^{6\times6}}$;
$\boldsymbol{u} \in { {\mathbb{R}}^{3}}$ is the open-loop control input; 
$\mathcal{J}(\cdot) \in {\mathbb{R}}_{+}$ is a generic cost functional that may account for time or control input minimization, among other factors;
$\boldsymbol{\Phi}(t_{i+1}, t_i) \in { {\mathbb{R}}^{6 \times 6}}$ denotes the state transition matrix of a natural dynamics from instant $t_i$ to $t_{i+1}$; and $\boldsymbol{B}(t_i) \in { {\mathbb{R}}^{6 \times 3}}$ is the control input matrix. 
A continuously differentiable function ${\mathcal{M}}^{*}(\boldsymbol{x}(t_i), \boldsymbol{X}(t_i)) \in {\mathbb{R}}$ represents the satisfaction of safety requirements on the uncontrolled trajectory after contingency, as a function of the controlled mean state and error covariance at contingency instant $(\boldsymbol{x}(t_i), \boldsymbol{X}(t_i))$ \cite{guffanti_jgcd_2023}.
The term $1-\Delta_{i}$ denotes the confidence level of constraint satisfaction at $t_i$.
The set of admissible control inputs is denoted as ${\mathcal{U}}$, while $\boldsymbol{x}_0 \in { {\mathbb{R}}^{6}}$ and $\boldsymbol{x}_f \in { {\mathbb{R}}^{6}}$ are the initial and terminal states, respectively. 

In this problem, the state covariance $\boldsymbol{X}(t_{i})$ is assumed to be a deterministic function of the state and control pair $f(\boldsymbol{x}(t_i), \boldsymbol{u}(t_i))$, emulating a closed-loop behavior where the spacecraft obtains a new state estimate with the error covariance at each time step. 
Details on the covariance model and its propagation are discussed in Sec. \ref{sec_IV}. 

As this optimization addresses open-loop control, the nominal (mean) dynamics, boundary conditions, and objective function are considered.
State uncertainty is incorporated into the path constraints, where the fault-tolerant behavior is expected to be achieved with an uncertainty margin. 

The probabilistic formulation of Eq. \ref{OCP} corresponds to an infinite-dimensional deterministic problem, which is generally intractable. 
To convert this into a finite-dimensional deterministic problem, the chance constraint is rewritten as follows \cite{guffanti_jgcd_2023}: 
\begin{equation}\label{PS_constr}
    {\mathcal{M}}^{*}(\boldsymbol{x}(t_i)) + \beta(\boldsymbol{X}(t_i), t_i) \leq 0.
\end{equation}
The term $\beta(\boldsymbol{X}(t_i), t_i) \in {{\mathbb{R}}_{\geq0}}$ is introduced into the evaluation of safety on the mean state ${\mathcal{M}}^{*}(\boldsymbol{x}(t_i)) := {\mathcal{M}}^{*}(\boldsymbol{x}(t_i), \boldsymbol{X}(t_i)=\boldsymbol{0})$ to account for the influence of state uncertainty on $\boldsymbol{x}(t_i)$ \cite{boyd_chance}. 
This decouples the contribution from the mean state and the uncertainty margin based on the state error covariance to evaluate the satisfaction of the chance constraint. 
Note that for the deterministic case, $\beta(\boldsymbol{X}(t_i)=\boldsymbol{0}, t_i)=0$.
By applying the transformation in Eq. \ref{PS_constr}, the problem presented in Eq. \ref{OCP} becomes a non-convex optimization problem.
One approach to solving this type of optimization problem is through a local direct optimization method such as sequential convex programming (SCP) \cite{malyuta_scp_2022}.

%%%%%%%%%%%%%%%%%%%%%%%%%%%%%%%%%%%%%%%%%%%%%
\section{Robust Trajectory Generation via Transformers}\label{sec_III}
%%%%%%%%%%%%%%%%%%%%%%%%%%%%%%%%%%%%%%%%%%%%%

This section first outlines relevant preliminaries of the Autonomous Rendezvous Transformer (ART) \cite{art_ieeeaero24} \cite{art_ral24}, followed by a presentation of the novel contributions in this work that enable the extension of the ART to a CC-OCP.

\subsection{Autonomous Rendezvous Transformer (ART)}

At its core, ART casts trajectory generation as a sequential prediction problem, whereby a Transformer model---trained on a dataset of previously collected near-optimal trajectories---returns a coherent sequence of states, controls, and performance parameters through autoregressive generation.
A key element for enabling reliable trajectory generation is the representation of trajectory data as sequences suitable for modeling by a Transformer.
Specifically, ART defines the following trajectory representation: 
\begin{equation}\label{eq:traj_representation}
    \tau = \{{\mathcal{P}}(t_1), \bx(t_1), \bu(t_1), \ldots,  {\mathcal{P}}(t_N), \bx(t_N), \bu(t_N) \},
\end{equation}
where ${\mathcal{P}}(t_i)$ denotes the vector of performance parameters at time $t_i$.
Generally, ${\mathcal{P}}(t_i)$ may include diverse metrics that define the desired performance of the trajectory.
In this work, two metrics are defined; ${{\mathcal{P}}(t_i)} = \left\{R(t_i), C(t_i) \right\}$, where $R(t_i) \in {\mathbb{R}}$ is the \textit{reward-to-go}, a metric representing future optimality, and $C(t_i) \in {\mathbb{N}}$ is the \textit{constraint-violation budget}, which indicates feasibility of the trajectory.
These metrics are defined as follows: 
\begin{align}\label{eq:rtg_ctg}
    R(t_i) & = -\sum_{j=i}^{N} {\mathcal{J}}(\boldsymbol{x}(t_j), \boldsymbol{u}(t_j), t_j), \\
    C(t_j) & = \sum_{j=i}^{N} {\textsf{C}}(t_j),
\end{align}
where ${\mathcal{J}}(\boldsymbol{x}(t_j), \boldsymbol{u}(t_j), t_j)$ represents the instantaneous reward that is directly associated with the cost function of the OCP (or CC-OCP), and ${\textsf{C}}(t_j) \in {\mathbb{N}}$ is the constraint violation metric incurred at instant $t_j$.

At a high level, three main processes are required to obtain a functional transformer model: dataset generation, model training, and test-time inference.

\textbf{Dataset Generation}: 
The first step in the proposed methodology involves generating a dataset suitable for effective Transformer training.
To accomplish this, a dataset of $N_d$ trajectories is generated by repeatedly solving diverse instances of the optimization problem in Eq. \ref{OCP}.
The raw trajectories are then rearranged according to the representation presented in Eq. \ref{eq:traj_representation}.
% An OCP (or CC-OCP) of interest is defined (Eq. \ref{OCP}). 
% A batch of this problem is solved using SCP with randomization of parameters (e.g., boundary conditions) to collect $N_d$ trajectories.
This approach facilitates the learning of a wide range of optimal trajectories along with their corresponding performance metrics. 
Specifically, exposing the Transformer to near-optimal trajectories with randomized initial conditions and performance metrics is crucial for enabling it to learn the effects of various input modalities. 
This diversity in training data allows the model to capture the underlying dynamics and variations in trajectory behavior, which is essential for generalizing to unseen scenarios and maintaining robust performance across different conditions.
It is important to note that information on uncertainty is not explicitly represented in the trajectory data but is implicitly embedded within the state, control histories, and the constraint-violation budget outlined below. 

\textbf{Training}: 
The transformer is trained by employing the standard teacher-forcing procedure commonly used in training sequence models.
Specifically, as in \cite{art_ieeeaero24}\cite{art_ral24} the following loss function is minimized:
\begin{equation}
\resizebox{0.95\linewidth}{!}{$\begin{aligned}
    {\mathcal{L}}(\tau) = \sum_{n=1}^{N_d} \sum_{i=1}^{N} \left( \Vert \boldsymbol{x}^{(n)}(t_i) - \hat{\boldsymbol{x}}^{(n)}(t_i) \Vert_2^2 + \Vert \boldsymbol{u}^{(n)}(t_i) - \hat{\boldsymbol{u}}^{(n)}(t_i) \Vert_2^2 \right),
\end{aligned}$}
\end{equation}
where $\hat{\boldsymbol{x}}$ and $\hat{\boldsymbol{u}}$ represent the predicted state and control generated by the Transformer, respectively.

\textbf{Inference}: 
Once trained, the Transformer model can be used to autoregressively generate an open-loop trajectory for warm-starting.
In particular, given an initial input in the form of an initial state and user-defined values for the performance metrics, trajectory generation is defined by the following procedure: 
(i) the Transformer predicts a control $\hat{\boldsymbol{u}}_1$, (ii) the state $\boldsymbol{x}_1$ and control $\hat{\boldsymbol{u}}_1$ are propagated to compute the next state $\boldsymbol{x}_2$, (iii) the performance parameters are updated, and (iv) the previous three steps are repeated until the entire trajectory has been generated.
In step (ii), the ``dynamics-in-the-loop'' procedure \cite{art_ieeeaero24} is employed, whereby the next state is computed using a dynamics model available onboard. 
This ensures that the generated trajectory is always dynamically feasible. 
The update of the performance parameters in step (iii) involves decreasing the reward-to-go $R(t_1)$ and constraint-violation budget $C(t_1)$, by the instantaneous reward ${\mathcal{J}}(\boldsymbol{x}(t_i), \boldsymbol{u}(t_i), t_i)$ and constraint violation ${\textsf{C}}(t_i)$.
Within this context, an effective strategy for selecting the performance parameters at initialization sets $R(t_1)$ as a (negative) quantifiable lower bound of the optimal cost and $C(t_1) = 0$, incentivizing the generation of near-optimal and feasible trajectories.
Ultimately, the generated trajectories serve as a warm-start for the SCP algorithm, which solves the non-convex trajectory optimization problem (Eq. \ref{OCP}), converging to a local optimum near the warm-start and ensuring constraint satisfaction. 

\subsection{Extension of ART to CC-OCP and Post hoc acceptance check}

This paper explores two novel aspects that stem from the above procedure. 
First, the constraint-violation budget is redefined to encapsulate chance constraints. 
Specifically, the constraint violation metric ${\textsf{C}}(\beta(\boldsymbol{X}(t_j),t_j))$ is defined as a binary indicator of the satisfaction of the chance constraint, using Eq. \ref{PS_constr}:
\begin{equation}\label{ctg2}
\resizebox{0.95\linewidth}{!}{$
    {\textsf{C}}(\beta(\boldsymbol{X}(t_j),t_j)) = 
\begin{cases}
    1 & \text{if } {\mathcal{M}}^{*}(\boldsymbol{x}(t_j)) + \beta(\boldsymbol{X}(t_j),t_j) \leq 0 \\
    0 & \text{otherwise}.
\end{cases}$}
\end{equation}
It should be noted that the state covariance history is not encoded in the trajectory representation in Eq. \ref{eq:traj_representation}.
However, the above constraint violation metric can implicitly inform that the mean states of the optimized trajectories of the CC-OCP have a probabilistic margin to path constraints.
This enables ART to understand the error covariance attached to the current state and predict the near-optimal control law that probabilistically satisfies the safety constraints. 

Secondly, a confidence metric is developed to evaluate the quality of the generated trajectory and determine whether it should be used as a warm-start to the SCP during autonomous operations. 
Specifically in safety-critical environments, it is crucial to ensure the reliability of the model's output and detect out-of-distribution scenarios. 
This motivates the development of criteria for accepting or rejecting the AI-generated output, which is embedded within the decision-making process.
Ultimately, the acceptance threshold should be tied to the final performance of the SCP's converged solution, allowing for the prediction of the SCP's outcome without the need to solve the non-convex optimization problem.

This work proposes a post hoc evaluation criterion based on the error between desired and realized performance parameters, namely:
\begin{equation}
    {\mathcal{E}}_{{\mathcal{P}}} = \left\{{\mathcal{E}}_{R}, {\mathcal{E}}_{C} \right\} = \left\{\dfrac{R_{\text{ART}} - R_{\text{des}}}{R_{\text{des}}}, C_{\text{ART}} - C_{\text{des}} \right\}.
\end{equation}

Concretely, this criterion assesses the difference between the performance of the generated trajectory and the desired performance that is used to condition the autoregressive generation process.
This quantitative evaluation provides a basis for deciding whether the trajectory meets the required standards for use in safety-critical operations.
If the post hoc check is successful, the ART-generated trajectory is used as the warm-start for the SCP; otherwise, a fallback trajectory is employed.

\begin{figure*}[ht]
    \centering
    \includegraphics[width=0.85\textwidth]{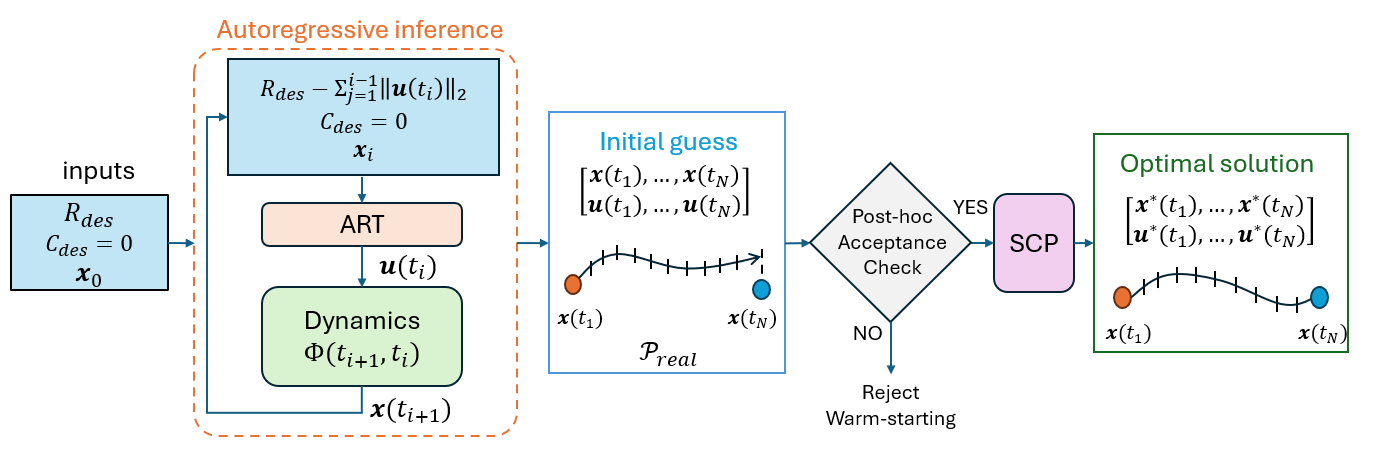}
    \caption{Pipeline of ART inference for non-convex chance-constrained trajectory optimization.}
    \label{fig:art_workflow}
\end{figure*}

Figure \ref{fig:art_workflow} summarizes an overview of the proposed pipeline.

%%%%%%%%%%%%%%%%%%%%%%%%%%%%%%%%%%%%%%%%%%%%%
\section{Rendezvous Scenario in Low Earth Orbit}\label{sec_IV}
%%%%%%%%%%%%%%%%%%%%%%%%%%%%%%%%%%%%%%%%%%%%%

\begin{table*}[ht] 
\centering
\caption{Simulation parameters for the rendezvous scenario in Low Earth Orbit.}\label{tab:scenario}
\resizebox{0.99\linewidth}{!}{
\begingroup
\renewcommand*{\arraystretch}{1.25}
\begin{tabular}{ c c c c c c c c c c c }
 \hline 
 \hline
  \multicolumn{2}{c}{Target orbit} & & \multicolumn{2}{c}{OCP} & & \multicolumn{2}{c}{SCP} \cite{oguri2023successive} & & \multicolumn{2}{c}{Dataset generation} \\
  \cline{1-2} \cline{4-5} \cline{7-8} \cline{10-11}
   $a$ [km] & 6738 & & horizon [orbits] & 1.29 & & [$\epsilon_\text{opt},\epsilon_\text{feas}$] & [1e-3,1e-3] & & $N_d$ & 115,000 \\
   $e$ [-] & 5.58e-4 & & $N_{\text{time}}$ [-] & 100 & & [$\alpha_1, \alpha_2$] & [2,3] & & Train split [\%] & 90 \\
   $i$ [deg.] & 51.64 & & $T_{\text{safety}}$ (orbits) & 1 & & $\beta_{\text{scvx*}}$  & 2 & & Test split [\%] & 10 \\
   $\Omega$ [deg.] & 301.04 & & $N_{\text{safety}}$ [-] & 20 & & $\gamma$ & 0.9 & &  & \\
   $\omega$ [deg.] & 26.18 & & $\boldsymbol{x}_{f, +T}$ & $\{ 0, 750 \textrm{m}, 0, 0, 0, 0\}$ & & [$r_{\text{min}}, r_{\text{max}}$] & [1e-6, 5e3]& & \multicolumn{1}{c}{Randomization ranges $\boldsymbol{x}_{0,rtn}$:} & (Uniform distribution) \\
   $M (t_1)$ [deg.] & 68.23 & & $\boldsymbol{x}_{f,rtn, -T}$ & $\{ 0, -750 \textrm{m}, 0, 0, 0, 0 \}$ & & $r^{(0)}$ & 100 & &  $\delta r_r$ (km) & $ \delta r_{r,0} + [-1.4, 1.4]$ \\
   period [h] & 1.528 & & $\boldsymbol{x}_{f,rtn, +R}$ & $\{ 750 \textrm{m}, 0, 0, 0, 0, 0 \}$ & & $w^{(0)}$ & 10 & & $\delta r_t$ (km) & $ \delta r_{r,0} + [-3, 3]$ \\
   & & & $\boldsymbol{x}_{f,rtn,-R}$ & $\{ -750 \textrm{m}, 0, 0, 0, 0, 0 \}$ & & max. iter. & 20 & & $\delta v_t$ (m/s) & $ \delta v_{t,0} + [-1, 1]$ \\
   & & & $\boldsymbol{x}_{0,rtn}$ & $\{ -4 \textrm{km} , -17.5 \textrm{km}, 0, 0, 6.849 \textrm{m/s}, 0 \}$ & & & & & & \\
   & & & $\{r_{r,\text{ae}},r_{t,\text{ae}},r_{n,\text{ae}}\}$ [m] & $ \{ 1000, 1000, 1000 \}$ & & & & & & \\
   & & & $\{r_{r,\text{koz}},r_{t,\text{koz}},r_{n,\text{koz}} \}$ [m] & $\{ 200, 200, 200 \}$ & & & & & & \\
   & & & $N_{\text{ae-to-koz}}$ & $55$ & & & & & & \\
  \hline
  \hline
\end{tabular}
\endgroup}
\end{table*}

The rendezvous scenario considered in this paper is a passively safe transfer of a binary system of spacecraft in LEO, from an initial to a target condition with four possible pre-docking waypoints ($\pm$ R-bar and $\pm$ V-bar) of the target. 
The scenario is inspired by \cite{berning2024} \cite{margolis2024robust}, which details are summarized in Table \ref{tab:scenario}. 
The fuel consumption is considered as a performance metric to be minimized in the optimization problem, where it is assumed that thrusters are always aligned with the direction of the desired control input. 
The reward-to-go is then defined as follows:
\begin{equation}
\resizebox{0.95\linewidth}{!}{$
\begin{aligned}
R(t_i) = - \sum_{j = i}^{N}{{\mathcal{J}}(\boldsymbol{x}(t_j), \boldsymbol{u}(t_j), t_j)} = -\sum_{j = i}^{N} || \Delta \boldsymbol{v}(t_j) ||_2.
\end{aligned}
$}
\end{equation}

\subsection{Relative Spacecraft Dynamics}

The relative motion between the servicer and the target is denoted by the state $\boldsymbol{x} \in {\mathbb{R}}^6$, and $\boldsymbol{u} \in {\mathbb{R}}^3$ is the control action that represents a maneuver applied by the actuation system.

An absolute reference orbit of the target spacecraft around a gravitational body in a two-body problem is uniquely defined by a set of orbital elements (OE): $\textbf{\oe} \in {\mathbb{R}}^6$. 
One representation is the quasi non-singular OE \cite{Vallado}: $\textbf{\oe} = \left\{ a, \nu, e_{x}, e_{y}, i, \Omega \right\}$, where $a$ is the semi-major axis, $\nu = M + \omega$ is the mean argument of latitude, $M$ is the mean anomaly, $\omega$ is the argument of periapsis, $\{e_x, e_y \} = \{ e \cos(\omega), e \sin(\omega) \}$ is the eccentricity vector, $e$ is the eccentricity, $i$ is the inclination, and $\Omega$ is the right ascension of the ascending node. 
The relative orbital motion of the servicer with respect to the target can be expressed equivalently using a relative Cartesian state, or using Relative Orbital Elements (ROE) which are nonlinear combinations of the OE of the servicer and the target \cite{damico_phd_2010}. 
The relative Cartesian state can be expressed in the Radial-Tangential-Normal (RTN) reference frame centered on the target that rotates around a primary body. 
The states are $\boldsymbol{x}_{rtn} = \{\delta \boldsymbol{r}_{rtn}, \delta \boldsymbol{v}_{rtn} \} \in {\mathbb{R}}^6$, where $\delta \boldsymbol{r}_{rtn} = \{\delta r_r, \delta r_t, \delta r_n \} \in {\mathbb{R}}^3$ and $\delta \boldsymbol{v}_{rtn} = \{\delta v_r, \delta v_t, \delta v_n \} \in {\mathbb{R}}^3$ denote the relative position and velocity in the RTN frame, respectively. 
The quasi-nonsingular ROE \cite{damico_phd_2010} is a nonlinear combination of the OE of the servicer (denoted with subscript $s$) and the target (denoted with no subscript) as: $\boldsymbol{x}_{roe} = \{ \delta a, \delta \lambda, \delta e_x, \delta e_y, \delta i_x, \delta i_y \} \in {\mathbb{R}}^6$, where $\delta a = (a_s - a)/a$ is the relative semi-major axis, $\delta \lambda = \nu_s - \nu + \left(\Omega_s - \Omega \right) \cos(i)$ is the relative mean longitude, $\{\delta e_x, \delta e_y \} = \delta e \{\cos(\varphi), \sin(\varphi) \} = \{e_s \cos(\omega_s) - e \cos(\omega), e_s \sin(\omega_s) - e \sin(\omega)\}$ is the relative eccentricity vector, and $\{\delta i_x, \delta i_y \}  = \delta i \{\cos(\phi), \sin(\phi) \} = \{i_s - i, \left(\Omega_s - \Omega \right) \sin(i)\}$ is the relative inclination vector. 
There exists both a nonlinear exact map and a first-order one-to-one mapping between relative Cartesian and ROE states. 
The latter is expressed as $\boldsymbol{x}_{rtn}(t) \approx \boldsymbol{\Psi}(t) \boldsymbol{x}_{roe}(t)$, where matrix $\boldsymbol{\Psi} \in {\mathbb{R}}^{6\times 6}$ is defined in \cite{damico_phd_2010} \cite{sullivan_nonlinear_2017} \cite{guffanti_jgcd_2023}.
For small spacecraft separation, the first-order linear map $\boldsymbol{\Psi}$ offers an accurate mapping of constraints between state representations.
The discrete-time linearized perturbed and controlled relative orbital dynamics of the servicer to the target is expressed more accurately in the ROE state compared to the RTN state, particularly for a large separation of spacecraft \cite{damico_phd_2010}. 
A state transition matrix $\Phi(t_{i+1}, t_i)$ that includes a variety of orbital perturbations is derived in \cite{koenig_2017_new} \cite{guffanti_linear_2019} \cite{guffanti_long-term_2017}, and the control input matrix $\boldsymbol{B}(t_i)$ is derived in \cite{gaias_imp_2015} \cite{chernick_2018_closed}, both as a function of the target's orbital elements. 

The control action is modeled as a net variation of the servicer's velocity (delta-v) in the RTN frame $\boldsymbol{u} = \Delta \boldsymbol{v}_{rtn} \in {\mathbb{R}}^3$.

In the rendezvous scenario, some constraints are more naturally formulated on the relative Cartesian state (e.g., terminal docking conditions and waypoints, safety-critical domains, etc.) \cite{pinglu_rpod}\cite{malyuta_rpod}\cite{malyuta_scp_2022}, whereas others are more advantageously formulated on the ROE state (e.g., perturbed orbital dynamics and passive safety conditions) \cite{koenig_2017_new}\cite{damico_proximity_2006}\cite{guffanti_jgcd_2023}. 
In the remainder of the paper, an ROE representation, which offers improved dynamics fidelity, is adopted for the formalization of the dynamics model in the OCP, i.e., $\bx = \bx_{roe}$.

\subsection{Uncertainty Modeling}

Three uncertainty sources are considered in this paper: navigation error, actuation error, and unmodeled system dynamics.
In particular, in this work, closed-loop behavior is modeled within open-loop optimization by assuming that the servicer acquires at each instant $t_i$ a new estimate of its relative state with a specified state error covariance. This uncertainty is then affected by the maneuver applied at $t_i$ (actuation error), and by the process noise (unmodeled system dynamics error) along the free-drifting trajectory after the application of the maneuver. 
These three sources of uncertainties are detailed in the following. 

\subsubsection{Navigation error} 
The filtered state estimate at $t_i$ is represented as a Gaussian distribution: ${\mathcal{N}}(\boldsymbol{x}^{\text{nav}}(t_i), \boldsymbol{X}^{\text{nav}}(t_i))$, where $\boldsymbol{x}^{\text{nav}}(t_i)$ and $\boldsymbol{X}^{\text{nav}}(t_i)$ denote the estimated mean and error covariance matrix, respectively.
It is assumed that at each iteration, the mean estimate is equivalent to the designed mean trajectory, i.e., $\boldsymbol{x}^{\text{nav}}(t_i) = \boldsymbol{x}(t_i)$.

The values of navigation uncertainty are based on the Absolute and Relative Trajectory Measurement System (ARTMS) \cite{kruger2024phd}, which is an angles-only navigation system that has been demonstrated on the NASA Starling  4x CubeSat swarm mission in LEO \cite{Kruger_STARFOX}. 
The navigation performance of ARTMS can be modeled as a function of the distance to the target as follows:
\begin{equation} \label{eq:artms_nav_general}
    \boldsymbol{X}^{\text{nav}}(\boldsymbol{x}(t_i)) = \frac{\rho}{a^2} \boldsymbol{s} \boldsymbol{s}^\top,
\end{equation}
where $\rho := ||\delta \boldsymbol{r}_{rtn}||_2$, and $\boldsymbol{s} \in {\mathbb{R}}_{\geq0}^{6}$ characterizes the navigation accuracy of each component of ROE.
In this work, $\boldsymbol{s}$ is modeled as a function of relative range as follows:
\begin{equation}\label{eq:artms_nav}
\resizebox{0.8\linewidth}{!}{$
\begin{aligned}
    \boldsymbol{s} = 
    \begin{cases}
        \boldsymbol{s}_{1}
        & \text{if } \rho > \rho_1,\\
        \dfrac{\boldsymbol{s}_{1} (\rho - \rho_2) + \boldsymbol{s}_{2} (\rho_1 - \rho)}{\rho_1-\rho_2} 
        & \text{if } \rho_2 < \rho < \rho_1, \\
        \boldsymbol{s}_{2} 
        & \text{otherwise},
    \end{cases}
\end{aligned}
$}
\end{equation}

where ($\rho_1, \rho_2$) and ($\boldsymbol{s}_{1}, \boldsymbol{s}_{2}$) are the pre-defined ranges and parameters.
When the relative range is in $\rho_2 < \rho < \rho_1$, the parameter is linearly interpolated. 

\subsubsection{Actuation error}
An actuation error of impulsive maneuvers is formulated based on the Gates model \cite{gates1963simplified}, which yields a zero-mean Gaussian noise in the frame $\{\boldsymbol{e}_1,\boldsymbol{e}_3,\boldsymbol{e}_3\}$ where $\boldsymbol{e}_1$ is the unit vector in direction of $\Delta \boldsymbol{v}$ and the $\boldsymbol{e}_2$ and $\boldsymbol{e}_2$ are in arbitrary-oriented orthogonal unit vectors. 
Specifically, the covariance $\boldsymbol{V} \in {\mathbb{R}}^{3 \times 3}$ is given as a function of $\Delta \boldsymbol{v}$ as follows: 
\begin{equation} \label{eq:gates_model}
\resizebox{0.95\linewidth}{!}{$
    \boldsymbol{V}(\Delta \boldsymbol{v}) = {\mathrm{Diag}}\{\sigma_r^2 + \Delta v^2 \sigma_s^2 , \sigma_a^2  + \Delta v^2 \sigma_p^2 , \sigma_a^2 + \Delta v^2 \sigma_p^2 \} ,
$}
\end{equation} 
where $\Delta v := || \Delta \boldsymbol{v} ||_2$, and $\{\sigma_s, \sigma_p,\sigma_r,\sigma_a\} \in {\mathbb{R}}_{\geq0}^{4}$ are model parameters.
Because it is assumed that $\bu = \Delta \boldsymbol{v}_{rtn}$, the covariance of the control input $\boldsymbol{U}\in {\mathbb{R}}^{3 \times 3}$ is given as 
\begin{equation} \label{eq:gates_model2}
% \resizebox{0.95\linewidth}{!}{$
    \boldsymbol{U}(\boldsymbol{u}(t_i)) = {\mathcal{R}}(\Delta \boldsymbol{v}(t_i)) \boldsymbol{V}(\Delta \boldsymbol{v}(t_i)) {\mathcal{R}}^\top(\Delta \boldsymbol{v}(t_i)),
% $}
\end{equation} 
where ${\mathcal{R}}(\Delta \boldsymbol{v}(t_i))$ is a rotation matrix from the coordinate frame $\{\boldsymbol{e}_1,\boldsymbol{e}_3,\boldsymbol{e}_3\}$ defined based on $\Delta \boldsymbol{v}(t_i)$ to the RTN frame. 

\subsubsection{Unmodeled system dynamics}
The discrepancy between the linearized dynamics model in Eq. \ref{OCP} and the true non-linear dynamics can be modeled as a discrete-time zero-mean white Gaussian process noise, where its covariance is denoted as $\boldsymbol{Q}(t_i) \in {\mathbb{R}}^{6 \times 6}$.
In this work, a time-invariant constant process noise matrix in the ROE state space is assumed.

\subsubsection{Linearized covariance dynamics}
With the above assumptions, the state error covariance $\boldsymbol{X}(t_{i})$ evolves in the following linearized dynamics:
%
% \small
\begin{equation}\label{DYND_Unc}
\resizebox{0.9\linewidth}{!}{$
\begin{aligned}
& \boldsymbol{X}(t_{i}) = \boldsymbol{X}^{\text{nav}}(\boldsymbol{x}(t_i)) + \boldsymbol{B}(t_{i}) \boldsymbol{U}(\boldsymbol{u}(t_i)) \boldsymbol{B}(t_{i})^\top + \boldsymbol{Q}(t_i) \\
& \boldsymbol{X}(t_{i+1}) = \boldsymbol{\Phi}(t_{i+1}, t_{i}) \boldsymbol{X}(t_{i}) \boldsymbol{\Phi}(t_{i+1}, t_{i})^\top.
\end{aligned}
$}
\end{equation}
% \small
%
It is important to note that the above equation shows that the state covariance history, $\boldsymbol{X}(t_{i})$, can be recovered from the nominal states and control inputs history,$(\boldsymbol{x}(t_i), \boldsymbol{u}(t_i))$.
Furthermore, since all uncertainties are modeled as Gaussian, the Gaussianity of the state covariance is preserved throughout the trajectory under the linearized dynamics model.

Details on the parameters considered in these uncertainty models can be found in Table \ref{tab:uncertainty}.

\begin{table*}[ht!]
\caption{Uncertainty Modeling}\label{tab:uncertainty}
\centering
\resizebox{0.55\linewidth}{!}{
\begingroup
\renewcommand*{\arraystretch}{1.25}
\begin{tabular}{l c l}
    \hline
    \hline 
    \centering
    Category  & Variables &  Value \\
    \hline
    Gates model \cite{gates1963simplified} & $\{\sigma_s, \sigma_p,\sigma_r,\sigma_a\}$ & \{2e-3, 3e-4,3e-4,3e-4\} \\ 
    \multirow{4}{*}{Navigation \cite{Kruger_STARFOX}} &
    $\boldsymbol{s}_{1}[\text{m}]$ &  [4e-5,4e-3,4e-5,2e-5,2e-5,4e-5]$^\top$  \\
    & $\boldsymbol{s}_{2} [\text{m}]$ &  [1e-4,4e-3,4e-5,2e-3,2e-3,2e-3]$^\top$ \\
    & $\rho_1 [\text{m}]$ & 10,000 \\
    & $\rho_2 [\text{m}]$ & 1,000 \\
   Unmodeled acc., & $\boldsymbol{Q}_{roe}, [\text{m}^2]$  & 1e-6 $\cdot \boldsymbol{I}_6 $  \\
    Confidence level & $q(\Delta_i)$ & 3.09 (cf. $\Delta_i = 0.001$) \\
    \hline 
    \hline 
\end{tabular}
\endgroup}
\end{table*}

\subsection{Safety Constraints with State Uncertainty}

Two safety constraints are considered in this rendezvous trajectory optimization problem: collision avoidance and passive safety.
Both constraints are evaluated based on the range between the two spacecraft, where the servicer is ensured to be excluded from safety domains around the target.  
In actual docking missions \cite{fehse2003automated}, multiple safety domains are employed, which progressively reduce in size as the servicer spacecraft approaches the target.
During the initial phase of rendezvous, the servicer's position is required to remain outside the Approach Ellipsoid (AE). 
As the relative range decreases, this safety boundary contracts into the Keep-Out Zone (KOZ) for the final approach phase. 
Both domains are modeled as an ellipsoid centered around the target, where the dimensions are characterized by a matrix $\boldsymbol{P} = {\mathrm{Diag}}\{1/r_r^2, 1/r_t^2, 1/r_n^2\} \in {\mathbb{R}}^{3 \times 3}$, where $\{r_r, r_t, r_n\}$ are the semi-major lengths in R/T/N directions, respectively. 
The epoch where the safety ellipsoid is switched from AE to KOZ is defined by a parameter $N_{\text{ae-to-koz}}$.
Table \ref{tab:uncertainty} lists the details of these values. 
 
As defined in Eq. \ref{OCP}, passive/free-drift safety is considered in this work, ensuring that an uncontrolled trajectory remains safe for at least $T_\text{safety}$ seconds after a contingency occurring at $t_i$. 
This safety requirement is enforced at $q$-$\sigma$ confidence level, applied to not only the state error covariance at $t_i$ but also the covariance propagated along the uncontrolled trajectory.
While a naive implementation of such a constraint would require numerous evaluations of inter-spacecraft separation along the uncontrolled arcs, the approach proposed in \cite{guffanti_jgcd_2023} allows for a single constraint evaluation to assess the safety of an uncontrolled trajectory originating from $t_i$. 
This evaluation is expressed as a function of the state and covariance at the contingency instant $(\boldsymbol{x}(t_i), \bar{\boldsymbol{X}}(t_i))$ in the form of Eq. \ref{PS_constr}. 
Due to the non-convexity of this constraint, sequential linearization around a reference trajectory ($\bar{\boldsymbol{x}}(t_i), \bar{\boldsymbol{u}}(t_i)$) and its corresponding covariance history $\bar{\boldsymbol{X}}(t_i)$ is introduced as follows \cite{malyuta_scp_2022}:
\begin{equation}\label{ps_lin}
    - \boldsymbol{a}^\top_{ps}(t_i) \boldsymbol{x}(t_i) + b_{ps}(t_i) + \beta(\bar{\boldsymbol{X}}(t_i),t_i) \leq 0,
\end{equation}
where
\begin{equation} \label{eq:chance_coeff}
\resizebox{0.95\linewidth}{!}{$\begin{aligned}
    &\boldsymbol{a}^\top_{ps}(t_i) = \bar{\boldsymbol{x}}^\top(t_i) \boldsymbol{\Phi}^\top(t^*,t_i) \boldsymbol{\Psi}^\top(t^*)\boldsymbol{D}^\top \boldsymbol{P} \boldsymbol{D} \boldsymbol{\Psi}(t^*) \boldsymbol{\Phi}( t^*,t_i) \in {\mathbb{R}}^{6} \\
    & b_{ps}(t_i) = \sqrt{\bar{\boldsymbol{x}}^\top(t_i) \boldsymbol{\Phi}^\top(t^*,t_i) \boldsymbol{\Psi}^\top(t^*)\boldsymbol{D}^\top \boldsymbol{P} \boldsymbol{D} \boldsymbol{\Psi}(t^*) \boldsymbol{\Phi}(t^*,t_i) \bar{\boldsymbol{x}}(t_i)} \in {\mathbb{R}}_{\geq0} \\
    & \beta(\bar{\boldsymbol{X}}(t_i),t_i) = q(\Delta_i) || \bar{\boldsymbol{X}}(t^*)^{1/2} \boldsymbol{a}_{ps}(t_i) ||_2 \in {\mathbb{R}}_{\geq0}.
\end{aligned}$}
\end{equation}
Here, $q(\Delta_i)$ refers to the percent point function (PPF) that indicates the desired confidence level. 
It is assumed that $\boldsymbol{X}(t_i)$ follows a Gaussian distribution (cf. Eq. \ref{DYND_Unc}), for which the PPF is readily available. 
The matrix $\boldsymbol{D} = [\boldsymbol{I}_3, \boldsymbol{0}_3] \in {\mathbb{R}}^{3 \times 6}$ is used to extract the relative position, $\delta \boldsymbol{r}_{rtn}$, from the state vector $\boldsymbol{x}_{rtn} = \{\delta \boldsymbol{r}_{rtn}, \delta \boldsymbol{v}_{rtn} \}$.
The instant $t^*$ represents the moment of the minimum distance between the target and the servicer in the free-drifting period $[t_i, t_i + T_\text{safety}]$, following the contingency. 
This evaluation is performed on the reference trajectory with the covariance state computed from it.
The covariance at the minimum separation instance, $\bar{\boldsymbol{X}}(t^*)$, is obtained by propagating $\bar{\boldsymbol{X}}(t_i)$ along the uncontrolled trajectory after contingency using Eq. \ref{DYND_Unc}. 
Notably, the collision avoidance constraint represents a special case of passive safety where $T_\text{safety} = 0$, a condition that lacks fault tolerance.  

\subsection{Sequential Convex Programming}

At each iteration of SCP, the following convex subproblem is solved: 
\begin{equation} \label{eq:cvx_subproblem}
\resizebox{0.95\linewidth}{!}{$
\begin{aligned}
& \text{given} & & \boldsymbol{a}_{ps}(t_i), b_{ps}(t_i), \beta(\bar{\boldsymbol{X}}(t_i),t_i) & \forall i \in [1, N] \\
& \underset{\boldsymbol{x}(t_i), \boldsymbol{u}(t_i), \xi_i}{\text{minimize}} & & \sum_{i = 1}^{N} || \boldsymbol{u}(t_i) ||_2  + \mathcal{J}_{\text{pen}}(\boldsymbol{\xi})&\\
& \text{subject to} & & \boldsymbol{x}(t_{i+1}) = \boldsymbol{\Phi}(t_{i+1}, t_{i}) (\boldsymbol{x}(t_{i}) + \boldsymbol{B}(t_{i}) \boldsymbol{u}(t_{i}) )  & \forall i \in [1, N) \\
&&& - \boldsymbol{a}^\top_{ps}(t_i) \boldsymbol{x}(t_i) + b_{ps}(t_i) + \beta(\bar{\boldsymbol{X}}(t_i),t_i) \leq \xi_i  & \forall i \in [1, N] \\
&&&  \boldsymbol{x}(t_1) =  \boldsymbol{x}_{0} \hspace{0.2cm} \textrm{and} \hspace{0.2cm} \boldsymbol{x}(t_N) =  \boldsymbol{x}_{f}, & \\
&&& \| \boldsymbol{x}(t_i) - \bar{\boldsymbol{x}}(t_i) \|_\infty \leq \rho &  \forall i \in [1, N] \\
&&& \| \boldsymbol{u}(t_i) - \bar{\boldsymbol{u}}(t_i) \|_\infty \leq \rho &  \forall i \in [1, N-1] \\
\end{aligned}$}
\end{equation}
where $\boldsymbol{\xi} = \{\xi_i\}_{i=1}^{N}$ represents a virtual buffer \cite{mao2018successive} that ensures the feasibility of the subproblem.
The utilization of the virtual buffer is further penalized in the objective function via $\mathcal{J}_{\text{pen}}(\boldsymbol{\xi})$.
No constraints on control input are imposed in this problem, assuming that the spacecraft possesses sufficient thrust capability for the scenario under consideration. 
Furthermore, the final two inequalities denote the (hard) trust region constraints, which ensures that each optimization variable are within the trust-region radius $\rho$ from the reference values. 

\begin{figure*}[h!]
    \centering
    \includegraphics[width=0.7\textwidth]{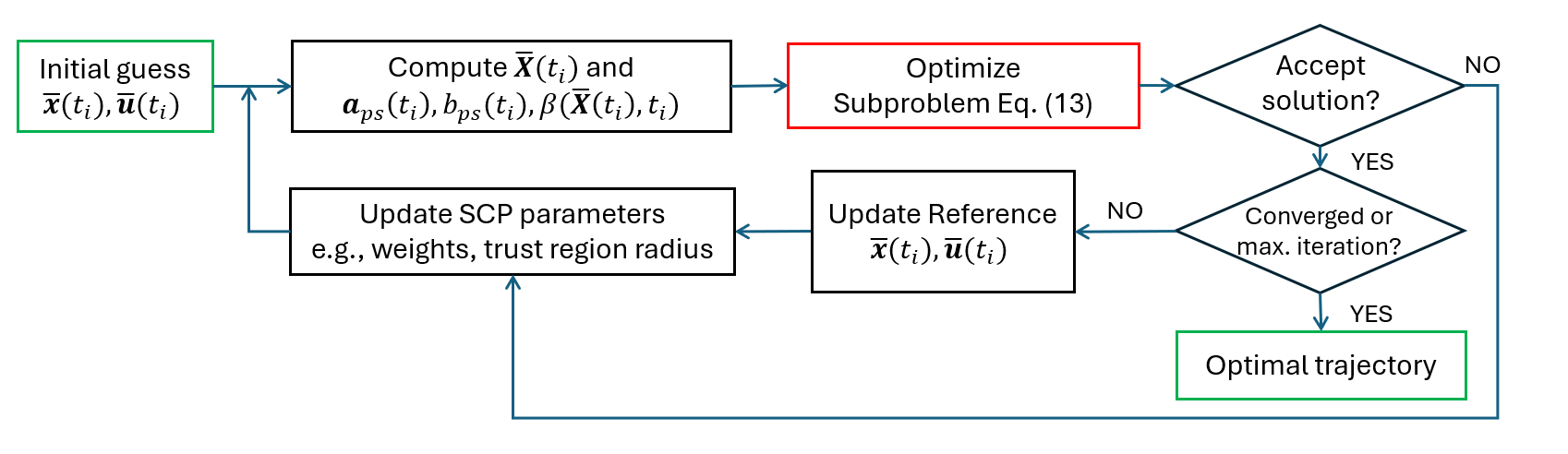}
    \caption{The workflow of SCP for CC-OCP. }
    \label{fig:scp_workflow}
\end{figure*}

This work adopts an SCP algorithm based on \texttt{SCVx*} \cite{oguri2023successive}.
For a detailed definition of the penalty term $\mathcal{J}_{\text{pen}}(\boldsymbol{\xi})$ and the update scheme of the parameters, the reader is referred to \cite{oguri2023successive}.
The workflow of the SCP is illustrated in Figure \ref{fig:scp_workflow}.
At each iteration, the reference state covariance history $\boldsymbol{X}(t_i)$ is recovered from the reference state and the control history $(\bx(t_i), \bu(t_i))$ using Eqs. \ref{eq:artms_nav_general}, \ref{eq:artms_nav}, \ref{eq:gates_model}, and \ref{DYND_Unc}.
Subsequently, Eq. \ref{eq:chance_coeff} provides the linearized chance constraint (cf. Eq. \ref{ps_lin}). 
Once the convex subproblem in Eq. \ref{eq:cvx_subproblem} is solved, the optimal solution is either accepted or rejected based on the evaluation of the nonlinear constraints. 
Finally, the reference and the SCP parameters (e.g., weights, trust region radius) are updated, if necessary. 
This iterative process continues until either the maximum number of iterations is reached, or convergence is achieved. 
%

%%%%%%%%%%%%%%%%%%%%%%%%%%%%%%%%%%%%%%%%%%%%%
\section{Experimental Results}\label{sec_V}
%%%%%%%%%%%%%%%%%%%%%%%%%%%%%%%%%%%%%%%%%%%%%
This section evaluates the performance of the proposed approach in challenging rendezvous scenarios.
The section begins with a detailed discussion of the dataset generation process and the transformer architecture employed.
Then, simulation results are presented, demonstrating statistically significant improvements in trajectory planning across diverse scenarios. 
Lastly, the impact of the post hoc acceptance check is assessed, along with qualitative analyses of selected trajectory examples.

\vspace{1mm}
\noindent\textbf{Dataset generation.} To generate a dataset of trajectories suitable for training the Transformer model, both the initial and terminal states are randomized within the domains specified in the last column of Table \ref{tab:scenario}.
Two datasets are generated: one where the trajectory optimization problem is solved in a deterministic fashion, and another where trajectory optimization is performed under chance constraints with provided uncertainty models.
Each dataset contains $N_d=115,000$ sample trajectories, which are divided into a $90$-$10\%$ train-test split.
Specifically, ART is trained on $90\%$ of the data, with the remaining $10\%$ reserved for evaluation, i.e., warm-start analysis.
% In this work, GPT model \cite{HuggFaceTransf} is used; readers are referred to the Appendix for the hyperparameters for the Transformer model. 

\begin{figure}
     \centering
     \begin{subfigure}[ht!]{\columnwidth}
         \centering
         \includegraphics[width=\columnwidth]{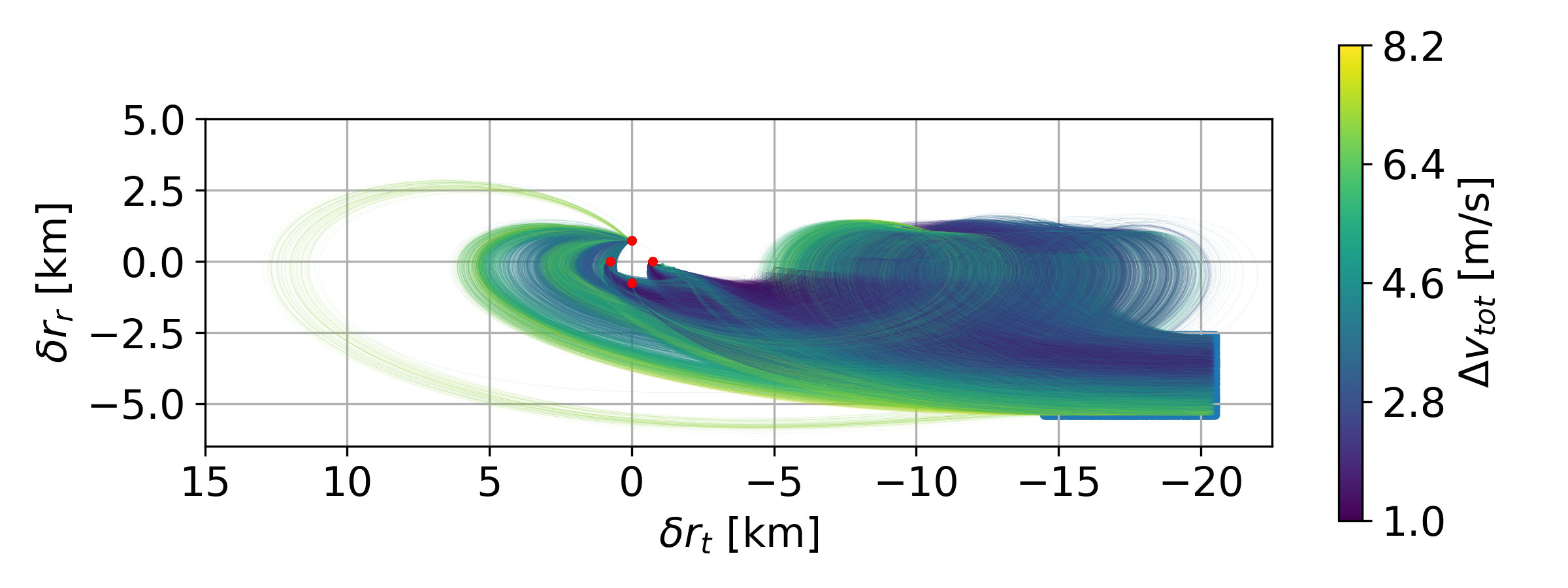}
         \caption{Deterministic trajectories}
         \label{fig:dataset_det}
     \end{subfigure}
     \\
     \begin{subfigure}[ht!]{\columnwidth}
         \centering
         \includegraphics[width=\columnwidth]{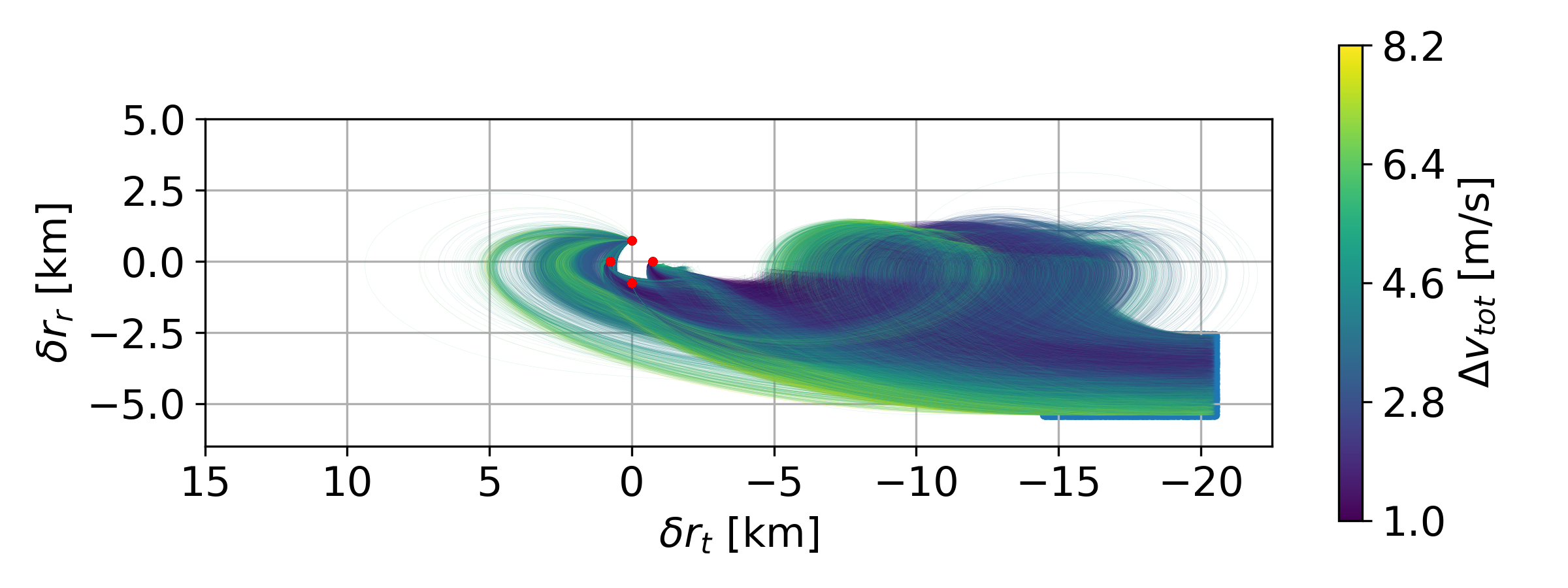}
         \caption{Chance-constrained trajectories}
         \label{fig:dataset_chance}
    \end{subfigure}
    \caption{Training datasets for the LEO rendezvous scenario ($N_d=115,000$).}
    \label{fig:dataset}
\end{figure}

The trajectories in each dataset are obtained through a structured sequence of optimization steps. 
First, a convex Two-Point Boundary Value Problem (TPBVP) is solved from the initial to the final condition, representing the simplest convex relaxation of the problem in Eq. \ref{OCP}.
Next, the convex solution is used to warm-start a non-convex problem that enforces a deterministic collision avoidance along the controlled trajectory (i.e., $\boldsymbol{X}(t_i)=\boldsymbol{0}$ and $T_\text{safety}=0$). 
Finally, the collision-avoidance solution is used to warm-start the problem in Eq. \ref{OCP}, which accounts for passive safety. 
The intermediate step of optimizing the deterministic collision avoidance trajectory is introduced due to the observed infeasibility of directly warm-starting the passive safety trajectory optimization with the convex TPBVP solutions, especially under chance-constrained scenarios (cf. Fig. \ref{fig:ws-analysis-chance}). 
This approach ensures the collection of a more diverse trajectory dataset, including trajectories not directly accessible from the convex TPBVP warm-start.
For each initial condition, the final feasible solution from the three trajectory optimizations is stored. 
Through this approach, each initial condition is associated with exactly one trajectory.
Figure \ref{fig:dataset} shows the trajectory datasets used in this paper in the radial/tangential (RT) plane, where the deterministic and chance-constrained trajectories are colored based on the total $\Delta v$ of the transfer (i.e., negative of the total reward-to-go), denoted as $\Delta v_{tot}$, in Figures \ref{fig:dataset_det} and \ref{fig:dataset_chance}, respectively. 

\vspace{1mm}
\noindent\textbf{ART architecture.} 
As in \cite{art_ieeeaero24}, the proposed model is a transformer architecture specifically designed to handle continuous inputs and outputs. 
Given an input sequence and a predefined maximum context length $K$, the model processes the last $5K$ sequence elements, corresponding to five modalities: target state, reward-to-go, constraint-violation budget, state, and control. 
Each modality’s sequence elements are projected through a linear transformation specific to that modality, producing a sequence of $5K$ embeddings.

Following the approach in~\cite{JannerEtAl2021}, an additional embedding is encoded for each timestep and added to the respective element embedding.
This sequence of embeddings is then fed into a causal GPT model consisting of six layers and six self-attention heads. 
The GPT architecture autoregressively generates a sequence of latent embeddings, which are projected through modality-specific decoders to predict the corresponding states and controls.

\vspace{1mm}
\noindent\textbf{Experimental design.} The warm-starting performance is evaluated by comparing the ART framework across several benchmarks, which are defined as follows:
\begin{itemize}
\item CVX: the solution of the convex TPBVP from initial to final state. This serves as the lower bound of the OCP and CC-OCP;
\item PS-CVX-D: the solution of the deterministic OCP in Eq. \ref{OCP} ($\boldsymbol{X}(t_i) = \boldsymbol{0}$) that is warm-started by CVX;
\item PS-ART-D: the solution of the deterministic OCP in Eq. \ref{OCP} ($\boldsymbol{X}(t_i) = \boldsymbol{0}$) that is warm-started by ART;
\item PS-CVX-C: the solution of CC-OCP in Eq. \ref{OCP} that is warm-started by CVX;
\item PS-ART-C: the solution of CC-OCP in Eq. \ref{OCP} that is warm-started by ART.
\end{itemize}

Note that although used in the dataset generation, collision avoidance trajectories are not used as warm-starts in the validation. 
This is because these solutions are not desirable onboard since they require the solution of an additional non-convex problem through SCP, which introduces further computational overhead.

\begin{figure*}[ht]
    \centering
    \includegraphics[width=\textwidth]{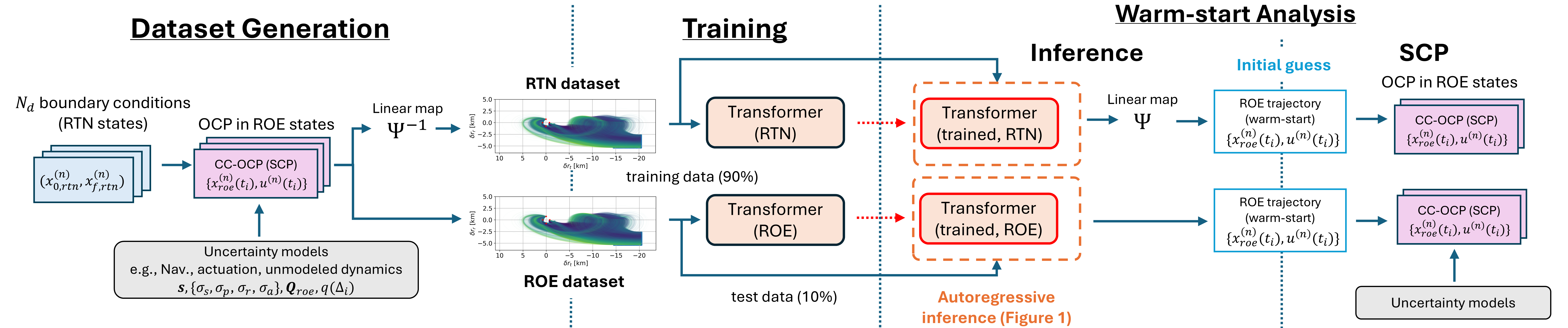}
    \caption{Pipeline of the experiment (dataset generation, training, inference, and warm-start analysis) with two state representations: RTN state (top) and ROE state (bottom). The two pipelines illustrate parallel approaches for training ART on either RTN or ROE state representations. This work assesses performance across both representations, allowing system designers to choose their preferred representation based on specific application needs.}
    \label{fig:entire_pipeline}
\end{figure*}

To evaluate the generalizability of ART to different state representations, two Transformer models are trained for each dataset.
Specifically, one model predicts the optimal control sequence based on the RTN state history, while the other generates trajectories from the ROE state history.
The workflow of this experiment, from dataset generation to training, inference, and warm-start analysis, is illustrated in Fig. \ref{fig:entire_pipeline} and outlined below.
First, the problem is defined with both initial and terminal states in the RTN coordinates (cf. Table \ref{tab:scenario}). 
The states are then converted into the ROE domain, where the CC-OCP is solved to generate the dataset, resulting in the ROE-based dataset that trains a Transformer model. 
In parallel, using the linear mapping $\Psi$, the optimized ROE trajectories are transformed into the RTN states, allowing for the training of the other Transformer model based on the RTN state history.
% solely trained using $\bx_{rtn}(t_i)$, along with the default model that is trained using $\bx_{roe}(t_i)$ that is directly obrained by the results of the dataset generation.
The inference performance of these two models is compared in the following subsection. \footnote{The experiments are conducted using AMD Ryzen 9 7950X CPU, and NVIDIA GeForce RTX 4090 GPU.}

% %%%%%%%%%%%%%%%%%%%%%%%%%%%%%%%%%%%%%%%%%%%%%
% \section{Autonomous Rendezvous Transformer Results} \label{sec_VI}
% %%%%%%%%%%%%%%%%%%%%%%%%%%%%%%%%%%%%%%%%%%%%%

\subsection{Warm-start Results}

\begin{figure}
     \centering
     \begin{subfigure}[ht!]{\columnwidth}
         \centering
         \includegraphics[width=\columnwidth]{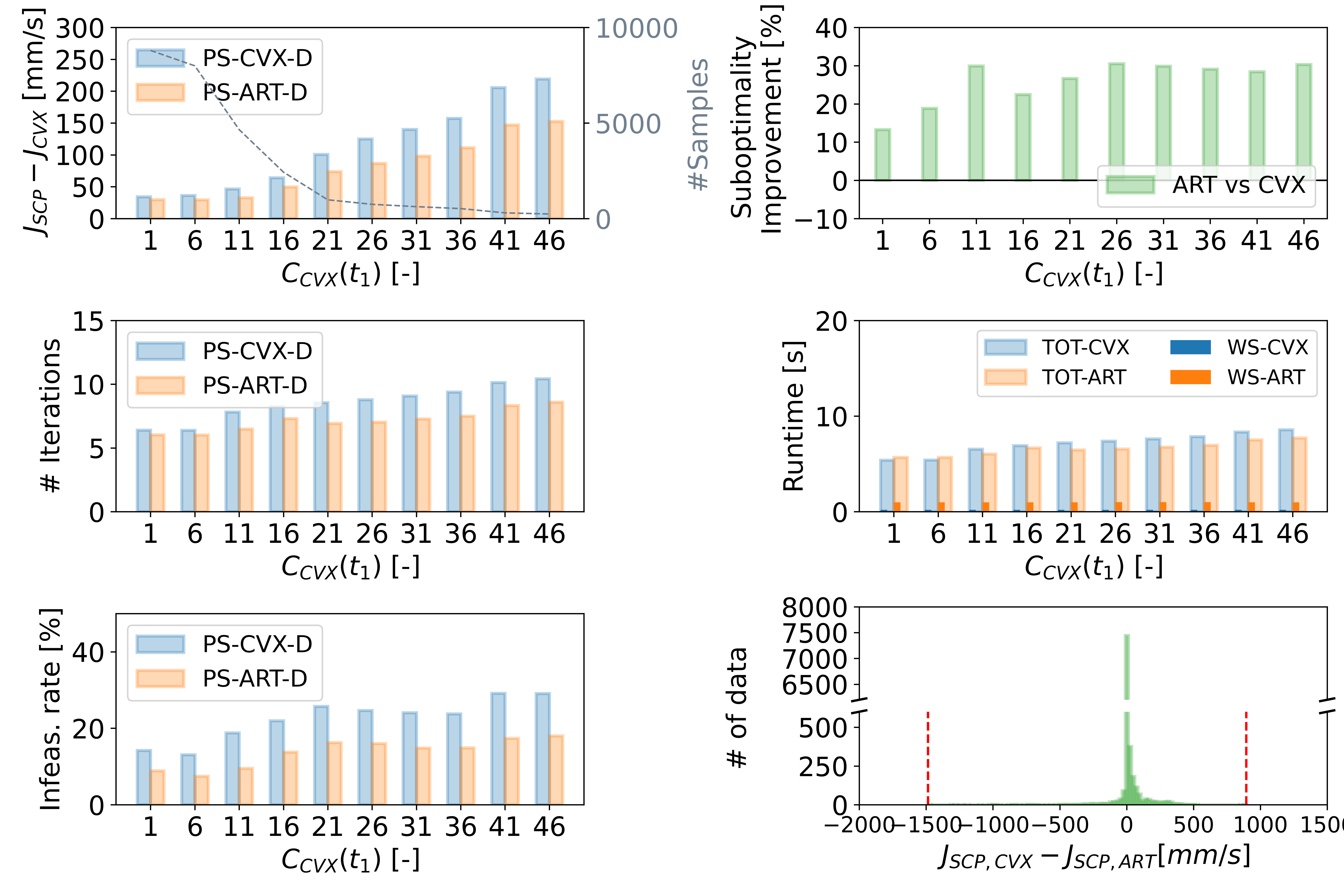}
         \caption{RTN}
         \label{fig:ws-analysis-det-rtn}
     \end{subfigure}
     \\
     \begin{subfigure}[ht!]{\columnwidth}
         \centering
         \includegraphics[width=\columnwidth]{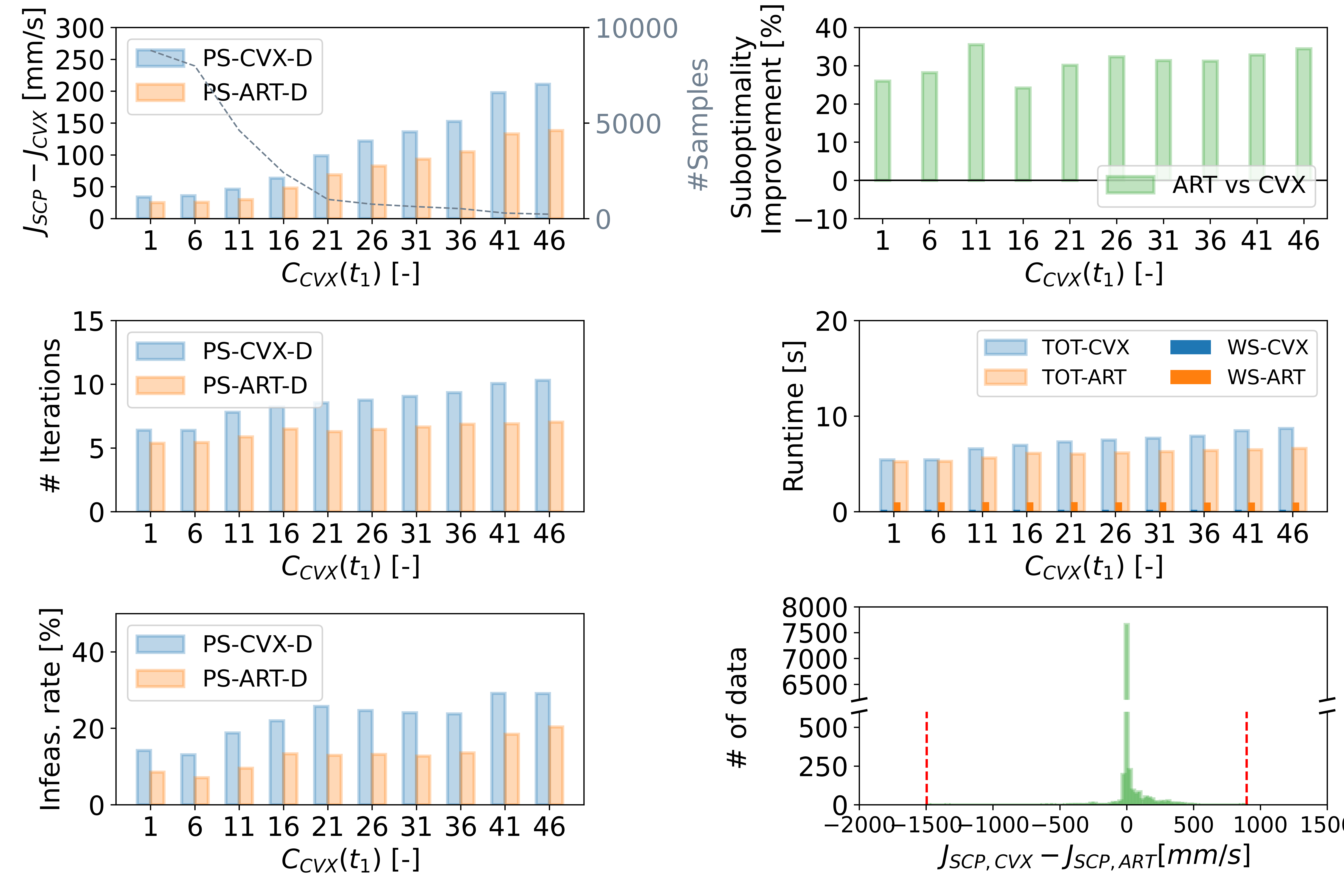}
         \caption{ROE}
         \label{fig:ws-analysis-det-roe}
    \end{subfigure}
    \caption{Performance comparison of SCP with different warm-starting strategies (deterministic scenario), plotted alongside the initial constraint-violation budget of CVX solutions, $C_\text{CVX}(t_1)$.}
    \label{fig:ws-analysis-det}
\end{figure}

\begin{figure}
     \centering
     \begin{subfigure}[ht!]{\columnwidth}
         \centering
         \includegraphics[width=\columnwidth]{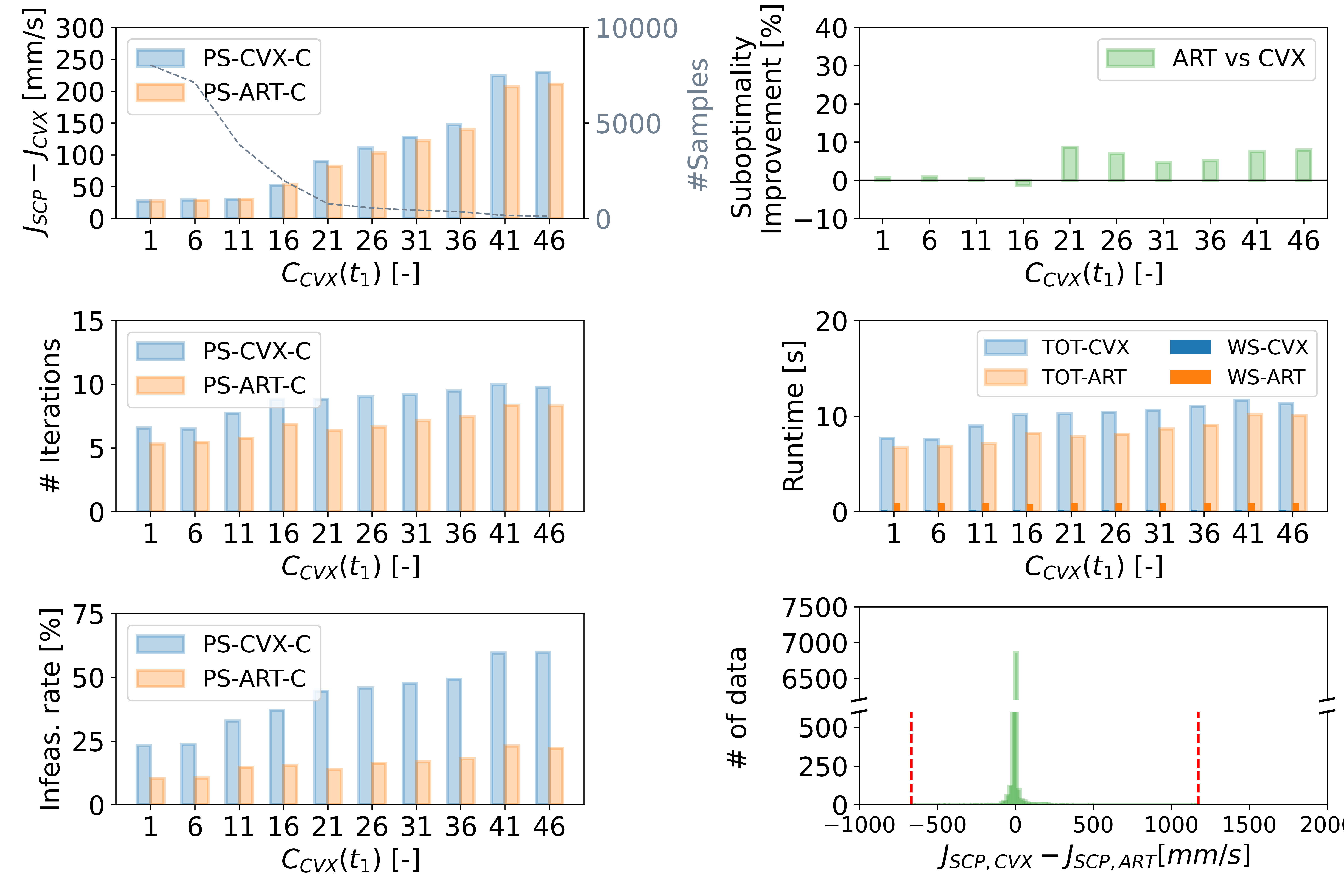}
         \caption{RTN}
         \label{fig:ws-analysis-chance-rtn}
     \end{subfigure}
     \\
     \begin{subfigure}[ht!]{\columnwidth}
         \centering
         \includegraphics[width=\columnwidth]{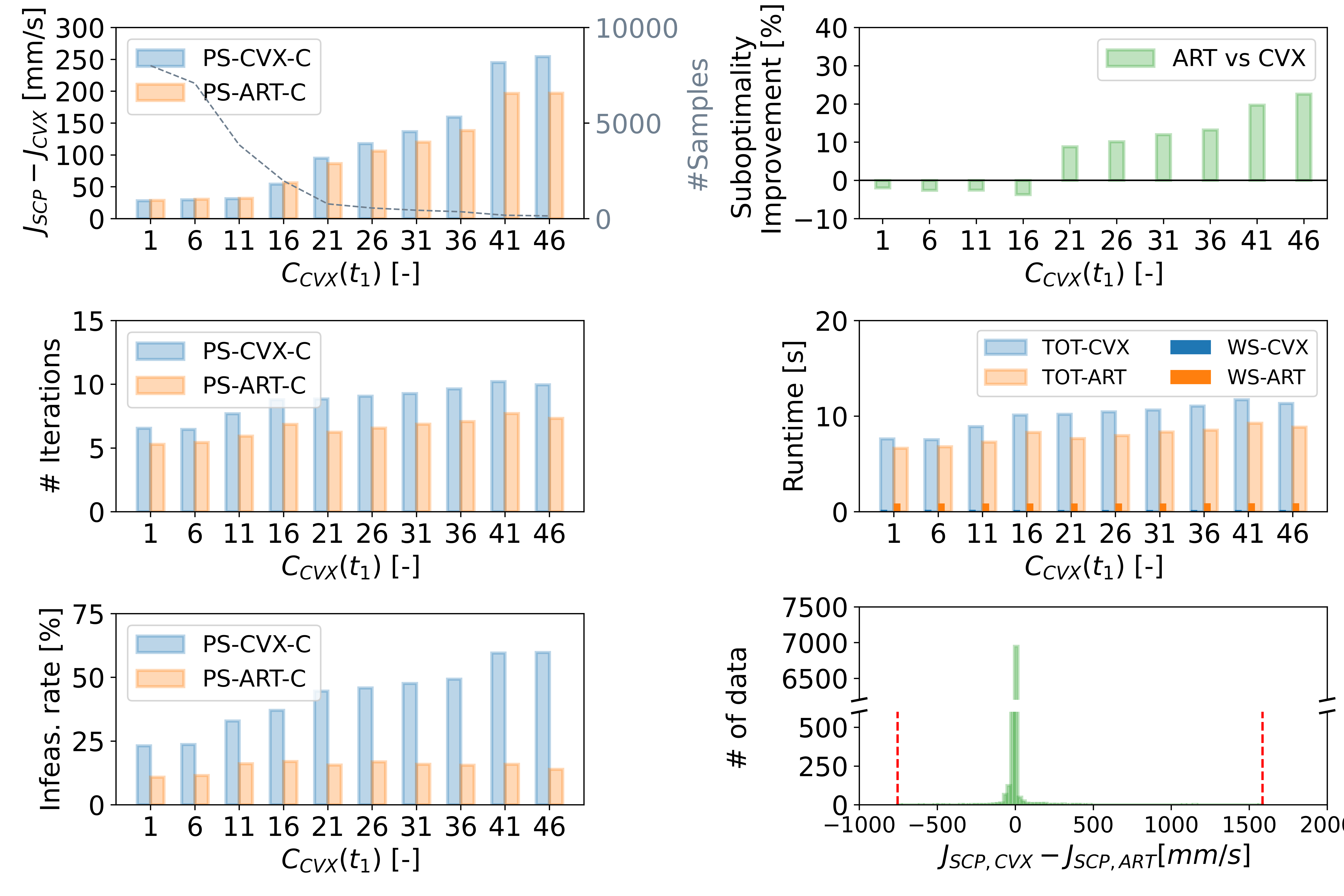}
         \caption{ROE}
         \label{fig:ws-analysis-chance-roe}
    \end{subfigure}
    \caption{Performance comparison of SCP with different warm-starting strategies (chance-constrained scenario), plotted alongside the initial constraint-violation budget of CVX solutions, $C_\text{CVX}(t_1)$.}
    \label{fig:ws-analysis-chance}
\end{figure}

Figures \ref{fig:ws-analysis-det} and \ref{fig:ws-analysis-chance} present the warm-starting performance for both the deterministic OCP and the robust CC-OCP, respectively. 
Figures \ref{fig:ws-analysis-det-rtn} and \ref{fig:ws-analysis-chance-rtn} correspond to the ART model trained on the RTN-based dataset, while Figures \ref{fig:ws-analysis-det-roe} and \ref{fig:ws-analysis-chance-roe} show the results for the ROE-based dataset. 
It is important to note that, in both cases, the OCP (CC-OCP) is ultimately solved using the ROE representation, as outlined in Fig. \ref{fig:entire_pipeline}. 
Each subplot compares key performance metrics of the SCP warm-started by ART (PS-ART) versus the CVX approach (PS-CVX), focusing on the following criteria: 
\begin{itemize}
    \item \textit{Absolute suboptimality improvement}: the difference between the converged objective values of the SCP and the CVX solution, i.e., $J_\text{SCP}-J_\text{CVX}$.
    \item \textit{Suboptimality improvement ratio}: calculated as $\frac{J_{\text{SCP,CVX}} - J_{\text{SCP,ART}}}{J_{\text{SCP,CVX}} - J_\text{CVX}}$, this ratio indicates the percentage improvement in suboptimality due to warm-starting with ART.
    \item \textit{Number of iterations required for SCP convergence}: evaluates the efficiency of each approach, with fewer iterations reflecting faster convergence.
    \item \textit{Runtime of the algorithm}: this includes the time required for both warm-starting and SCP solution.
    \item \textit{Infeasibility rate}: the percentage of solutions that were either infeasible or failed to converge.
\end{itemize}

The majority of the subplots in Figures \ref{fig:ws-analysis-det}, \ref{fig:ws-analysis-chance} are plotted as a function of the initial constraint-violation budget of CVX solutions, $C_\text{CVX}(t_1)$, along the x-axis.
Specifically, each point on the x-axis represents an average of the results for all $C_\text{CVX}(t_1)$ greater or equal to that value. 
Intuitively, higher values of $C_\text{CVX}(t_1)$ correspond to more challenging scenarios, whereby the CVX solution exhibits more frequent constraint violations.
In such cases, the optimization problem becomes more complex, requiring greater adjustments to bring the solution within feasible bounds. 

Lastly, the bottom right subplot displays the histogram of the optimality gap between the converged values of PS-CVX and PS-ART (i.e., $J_{\text{SCP,CVX}} - J_{\text{SCP,ART}}$) across the test trajectories. 
In this plot, a positive value indicates that the PS-ART has achieved a more optimal solution compared to the PS-CVX. 
Additionally, two red dotted lines mark the maximum and minimum values of the optimality gap across the entire test set. 

Concretely, the improvement in suboptimality achieved through warm-starting with ART is consistently demonstrated in both deterministic and chance-constrained scenarios.
As observed in the top rows of Figures \ref{fig:ws-analysis-det}, \ref{fig:ws-analysis-chance}, these results are confirmed across both RTN and ROE representations for training and inference.
The benefits are particularly pronounced in trajectories with higher values of $C_\text{CVX}(t_1)$, achieving up to 35\% and 22\% cost improvements for the deterministic and chance-constrained scenarios, respectively.
When $C_\text{CVX}(t_1)$ is low, it indicates that the CVX is close to being passively safe, thereby reducing the comparative advantages of using ART for simpler problems, as the converged solutions of PS-CVX and PS-ART tend to be similar.  
Conversely, significant improvements are observed when the ART warm-start is applied to highly infeasible CVX solutions, highlighting the effectiveness of ART in addressing more challenging trajectory optimization problems.

In the middle rows of Figures \ref{fig:ws-analysis-det}, \ref{fig:ws-analysis-chance}, the results further indicate that warm-starts from ART converge to local optima more rapidly than those from CVX, resulting in fewer SCP iterations across the range of $C_\text{CVX}(t_1)$.
The middle right subplots show that inference using the Transformer model (WS-ART) requires more computational time than solving a single convex TPBVP (CVX). 
However, the reduced number of iterations in the subsequent SCP results in a comparable or shorter overall runtime, highlighting the advantageous performance of ART.

Additionally, the lower-left plots further indicate that the infeasibility rate improves with the ART warm-start.
In particular, the initial guess provided by the ART effectively mitigates the increase in the infeasibility rate in more challenging scenarios (i.e., high $C_\text{CVX}(t_1)$) domain. 
It is observable that the warm-starting based on the CVX solution suffers from a high infeasibility rate in the chance-constrained case compared to the deterministic case. However, the infeasibility rate of the ART warm-start stays relatively constant, presenting a compelling case for adopting a learning-based warm-starting strategy.  

Lastly, the lower right subplots show how most test trajectories result in an optimality gap of the converged solutions (i.e., $|J_\text{SCP,CVX} - J_\text{SCP,ART}|$) less than 0.25m/s.
However, there are instances where the ART warm-start leads to both significant improvements and deteriorations in optimality, with some cases exhibiting differences greater than 1 m/s.
% This variability is attributed to the mechanics of the SCP algorithm, it shows the importance of having a further post hoc acceptance check before deploying the control obtained by the SCP.

Notably, similar results emerge when comparing the results obtained using different state representations (RTN and ROE) for training the model, suggesting that ART can effectively generate high-performance trajectories across different state representations. 
In practice, the results using the ROE state exhibit slightly superior performance than those using the RTN state across the metrics considered. 

\subsection{Qualitative Assessment}

\begin{figure}[ht!]
    \centering
    \includegraphics[width=\columnwidth]{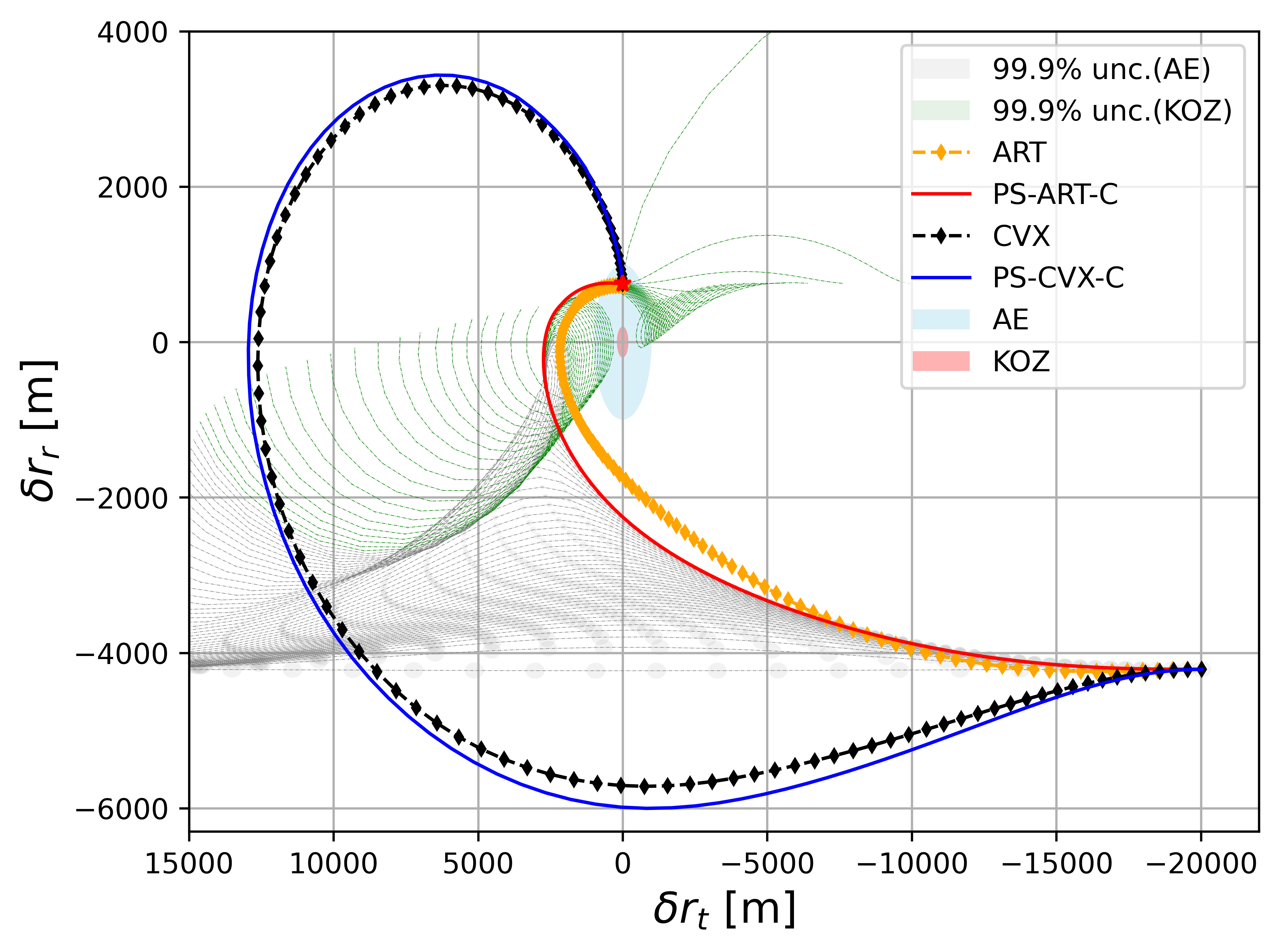}
    \caption{Sample chance-constrained trajectories with corresponding warm-starts. Green and gray trajectories are the free-drift trajectories from PS-ART-C with state covariance attached, which shows the passively safe behavior to AE and KOZ, respectively. }
    \label{fig:samlple_traj}
\end{figure}

Figure \ref{fig:samlple_traj} illustrates a representative trajectory that achieves one of the largest gains in optimality through the ART warm-start.
The figure presents the initial guesses of CVX, ART, and the locally optimal trajectories derived from the corresponding warm-start. 
The green and gray trajectories are propagated along the natural dynamics from the controlled trajectory (PS-ART-C) at each time step, ensuring passive safety for the AE and KOZ, respectively. 
Passive safety is achieved with the uncertainty margin provided by the uncertainty ellipsoids attached to each drifting trajectory.
First, the CVX solution has the optimal cost of $\Delta v_{tot} = 4.7916$ m/s. 
However, enforcing passive safety (PS-CVX-C) incurs an additional cost in the subsequent SCP convergence, leading to $\Delta v_{tot} = 5.2892$ m/s. 
On the contrary, the trained Transformer model returns a passively safe solution with $\Delta v_{tot} = 4.7895$ m/s, which accurately predicts the lower bound of the optimal cost (i.e., CVX).
Although this solution does not satisfy the terminal condition, this initial guess is further optimized by the following SCP, resulting in the total cost of the converged solution (PS-ART-C) of $\Delta v_{tot} = 4.7917$ m/s, almost matching the lower bound of the optimal cost derived from CVX.
Obtaining an initial guess comparable to the ART solution--namely, one that attains a low cost while satisfying the nonconvex constraint—-is substantially more challenging when relying on a relaxed convex formulation. This is because the relaxed convex problem does not incorporate information on the sensitivity of the additional cost imposed by enforcing nonconvex constraints. Consequently, this representative case highlights the advantage of Transformer-based initial guess generation, which accurately predicts transfer cost and safety requirements, thus yielding a high-quality initial guess in a highly multimodal solution space.

\begin{figure}[ht!]
    \centering
    \includegraphics[width=\columnwidth]{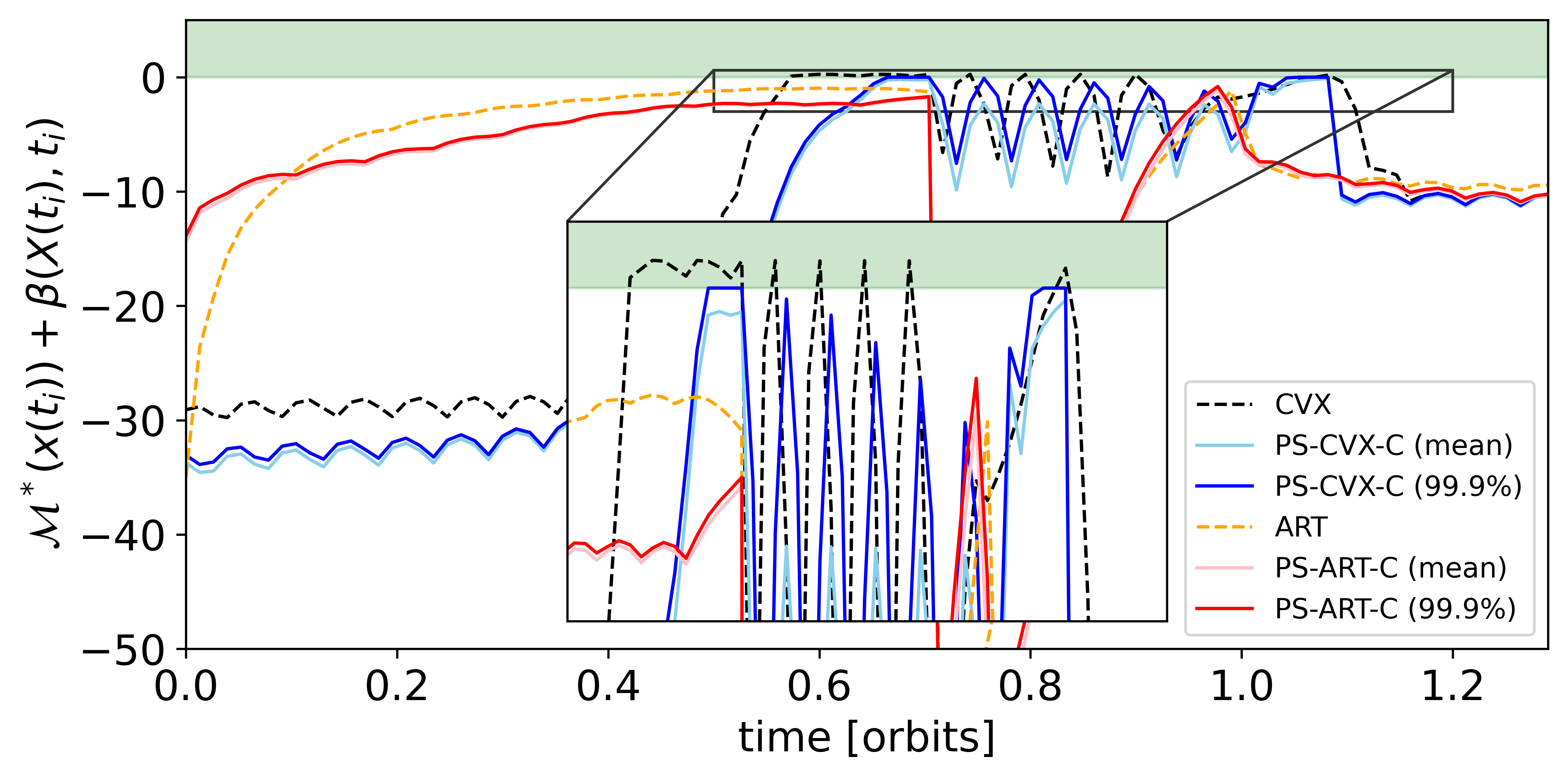}
    \caption{History of passive safety constraint satisfaction. The green zone represents the constraint violation.}
    \label{fig:samlple_traj_constr}
\end{figure}

Furthermore, Figure \ref{fig:samlple_traj_constr} illustrates the constraint satisfaction of the trajectories with a 99.9\% uncertainty margin for both PS-CVX-C and PS-ART-C.
The history of the left-hand side of Eq. \ref{ps_lin} is shown, where light colors represent the mean constraint values (i.e.,${\mathcal{M}}^{*}(\boldsymbol{x}(t_i)$) and dark colors represent the constraint value that accounts for the probabilistic margin. 
As defined in Eq. \ref{ps_lin}, the non-positive domain corresponds to constraint satisfaction. 
The plot illustrates that constraint violation in the warm-start trajectory from CVX is resolved in both locally optimal trajectories (PS-CVX-C), even when accounting for the 3-$\sigma$ margin of uncertainty. 
Furthermore, the generated ART warm-start is passively safe, satisfying the nonconvex constraint with the uncertainty margin before the solution is further optimized by the following SCP (PS-ART-C).

\subsection{Post hoc Acceptance Check on Transformer Output}

\begin{figure}
     \centering
     \begin{subfigure}[ht!]{\columnwidth}
         \centering
         \includegraphics[width=\columnwidth]{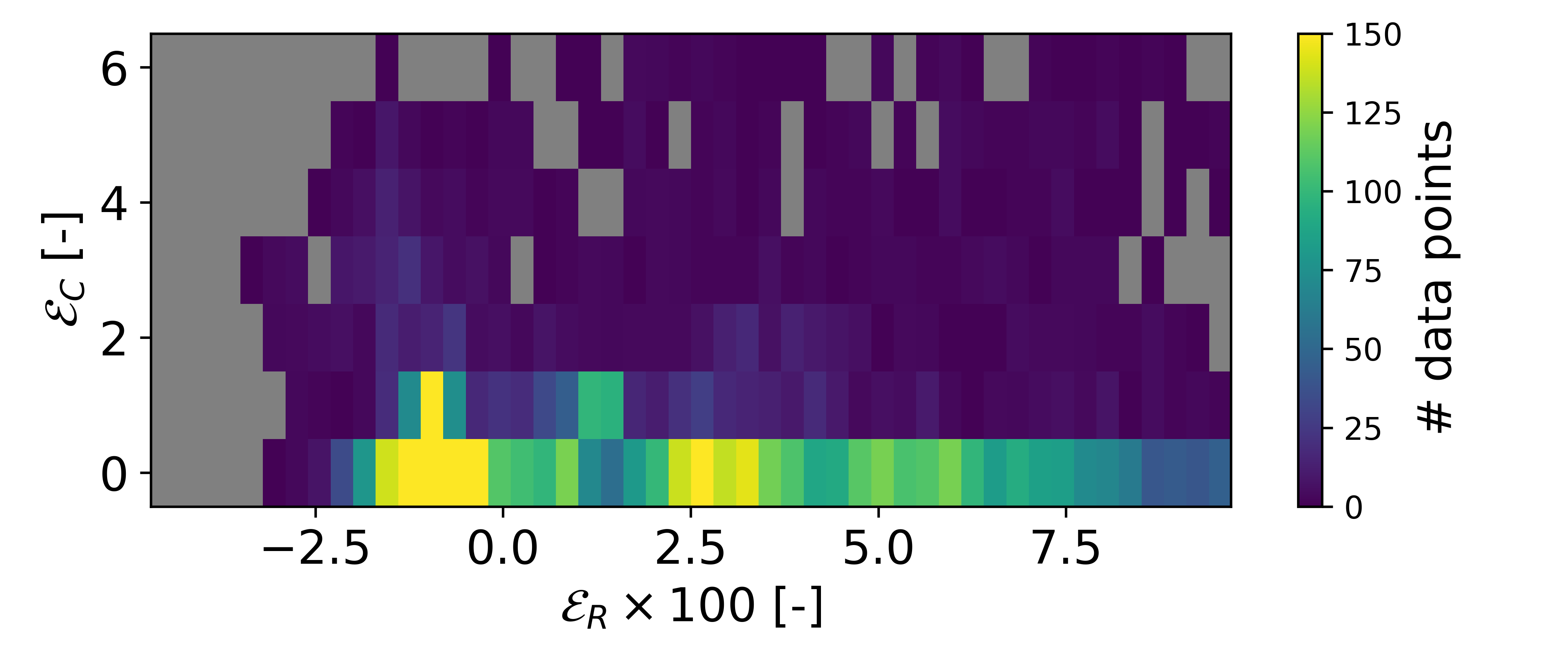}
         \caption{Data population}
         \label{fig:dataset_pop}
     \end{subfigure}
     \\
     \begin{subfigure}[ht!]{\columnwidth}
         \centering
         \includegraphics[width=\columnwidth]{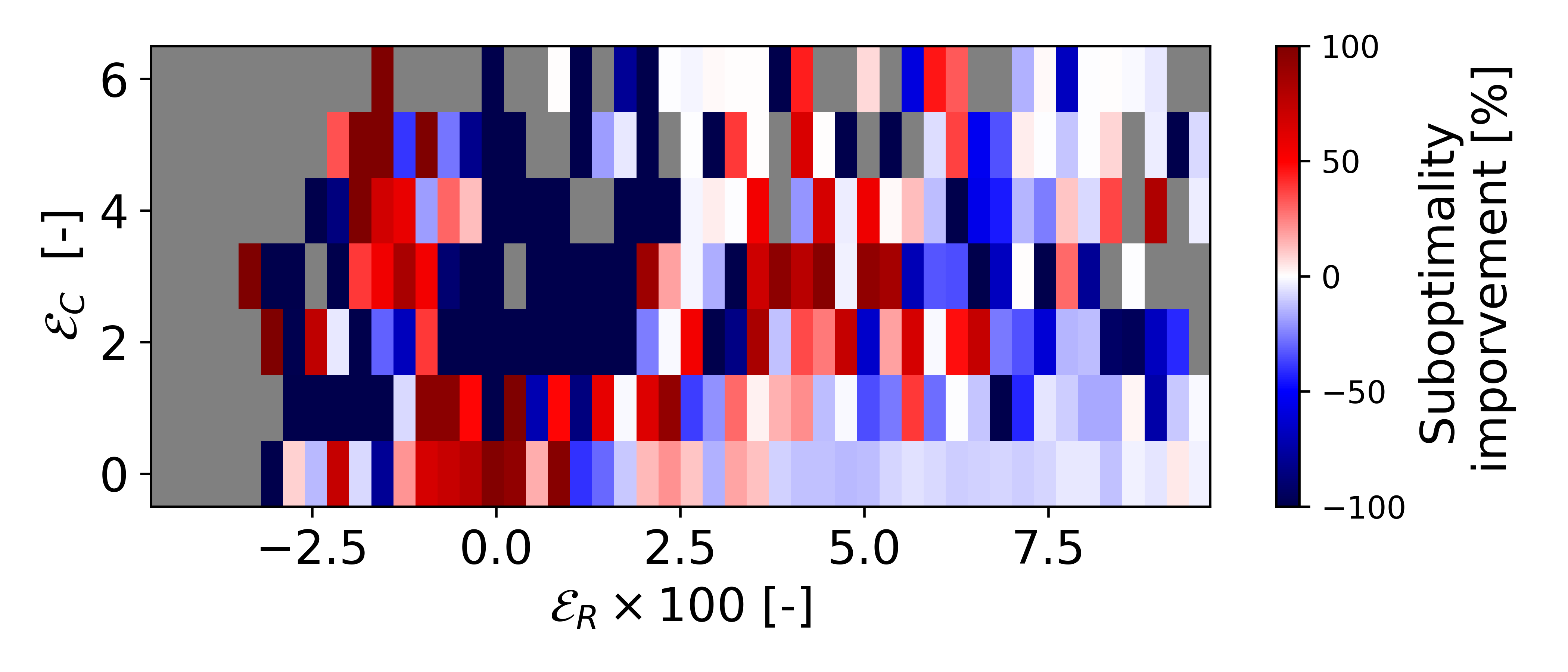}
         \caption{Suboptimality improvement}
         \label{fig:dataset_subopt_imprv}
    \end{subfigure}
    \\
     \begin{subfigure}[ht!]{\columnwidth}
         \centering
         \includegraphics[width=\columnwidth]{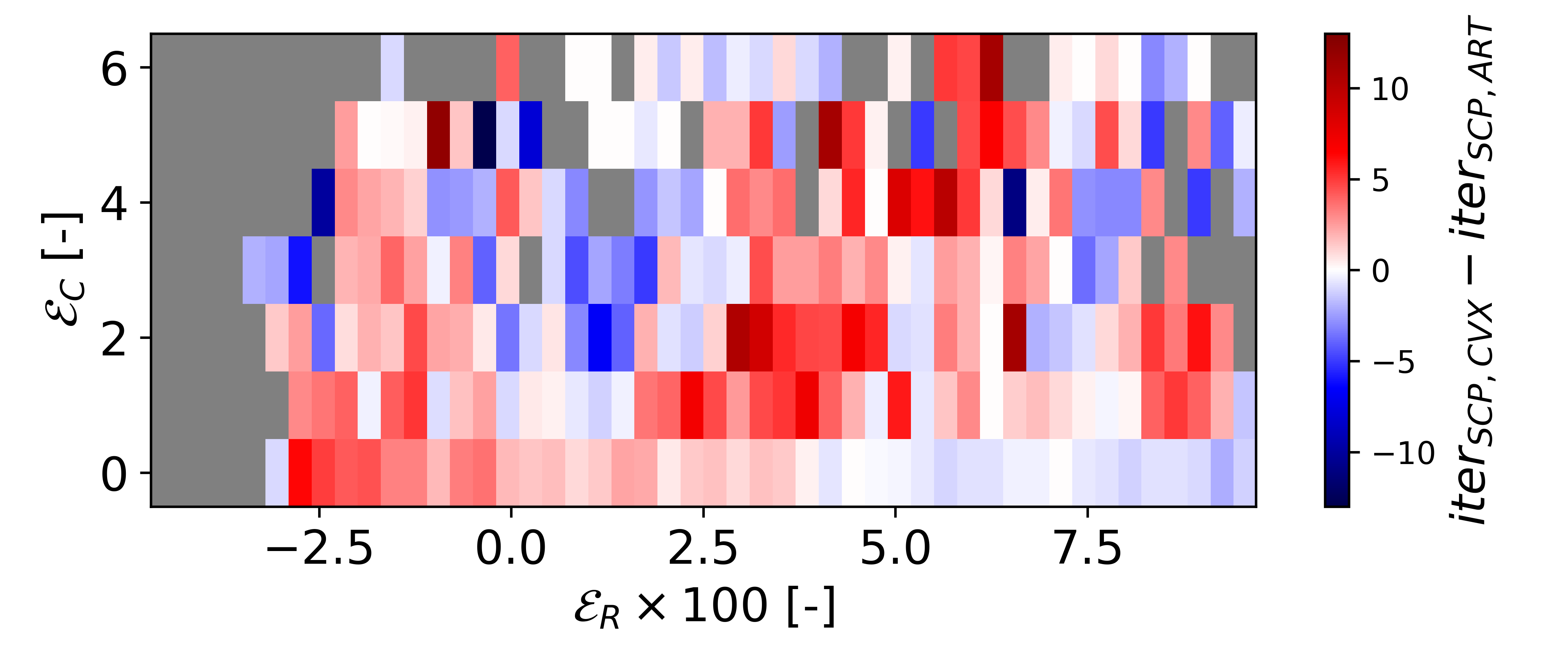}
         \caption{Iteration improvement}
         \label{fig:dataset_scp_iter}
    \end{subfigure}
    
    \caption{Analysis on the ART-generated trajectories in the warm-start analysis (chance-constrained, RTN-based dataset)}
    \label{fig:dataset_property}
\end{figure}

Finally, the post hoc acceptance check is discussed to filter out low-quality trajectories that should not be used as warm-starts for the SCP. 
This subsection focuses on a representative case from the chance-constrained RTN-based dataset.
First, Figure \ref{fig:dataset_pop} illustrates the population distributions of the trajectories used in the warm-start analysis.
Figures \ref{fig:dataset_subopt_imprv} and \ref{fig:dataset_scp_iter} depict the suboptimality improvement and reduction in the number of iterations, respectively, within the space defined by the post hoc criterion $\{ {\mathcal{E}}_R, {\mathcal{E}}_C \}$. 
In both subfigures, the red domain highlights the instances where the ART warm-starting outperforms the CVX warm-starting. 
Gray tiles indicate regions without data points. 

Several observations emerge from these figures. 
First, it is observable that the majority of the generated trajectory achieves passive safety, i.e., ${\mathcal{E}}_C = C_\text{ART}(t_1) = 0$. 
Moreover, given the current input definition (where the lower bound, i.e., CVX solution, serves as the desired reward-to-go), most generated trajectories result in costs below the objective of the original CVX solutions.  
However, it is important to note that these trajectories are generally infeasible, as they do not satisfy the terminal state condition. 
Additionally, in subfigure (b), it can observed that the suboptimality improvement ratio decreases as ${\mathcal{E}}_C$ and ${\mathcal{E}}_R$ increase. 
As expected, when the generated trajectory is (almost) passively safe and returns a near-optimal solution, there is a greater likelihood of the subsequent SCP converging to a more globally optimal extremum. 
Finally, the last subfigure illustrates the distribution of ART-generated trajectories which tend to converge faster than those initialized with CVX. 
Specifically, for trajectories demonstrating near-optimal performance, more iterations are often required to reach local suboptimality when the initial constraint-violation budget is high.
% although this area tends to outperform in terms of the suboptimality improvement by ART warm-start.
In contrast, trajectories that have high values of ${\mathcal{E}}_R$ and ${\mathcal{E}}_C$ (e.g., $2.5 < {\mathcal{E}}_R \times 100 \leq 7.5, 1 < {\mathcal{E}}_C \leq 4$) tend to reduce the number of iterations required for SCP convergence, even though the converged solutions may not necessarily be close to global optima.

\begin{figure}[ht!]
    \centering
    \includegraphics[width=\columnwidth]{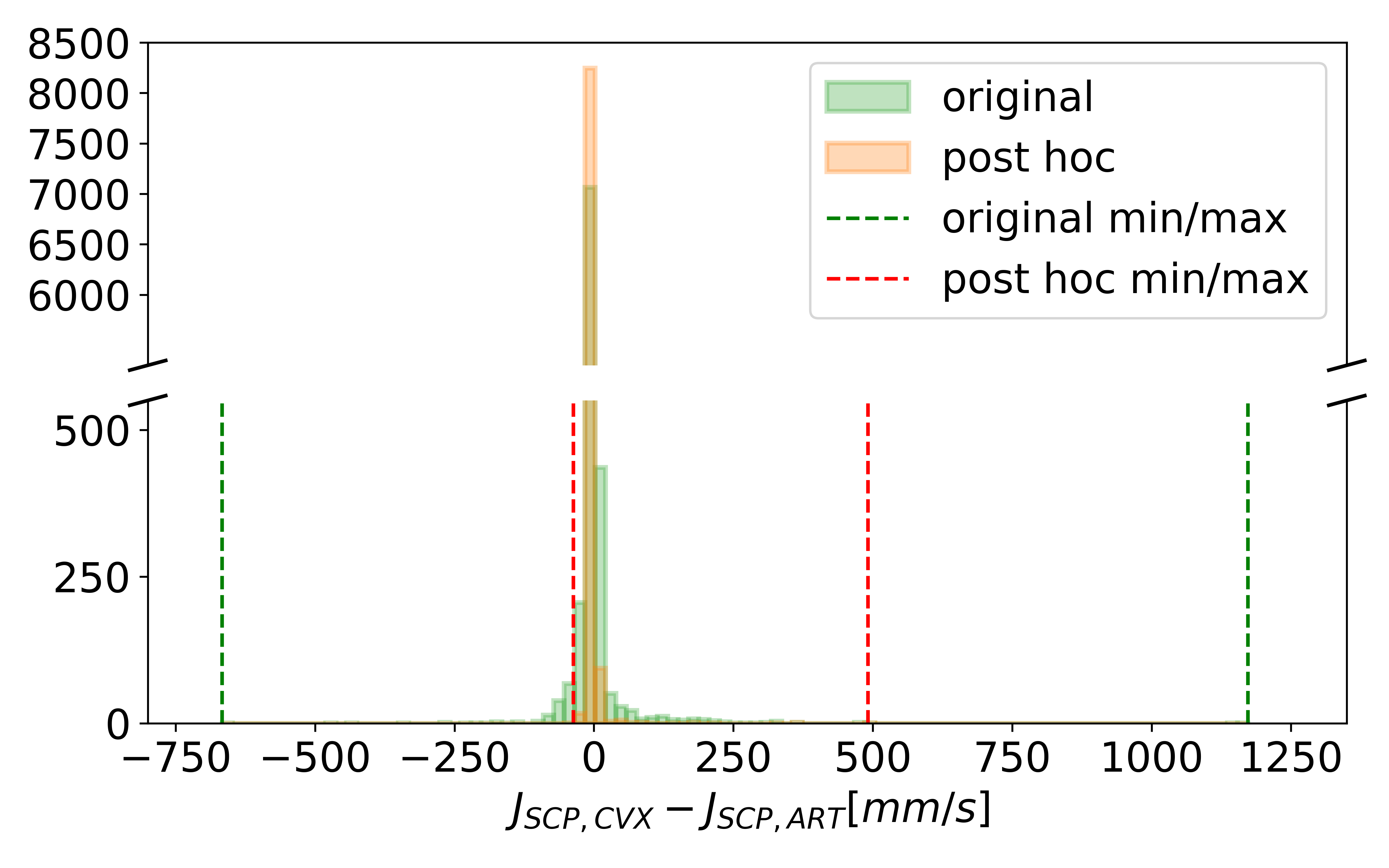}
    \caption{Hisogram of the solutions of the warm-start analysis with post hoc acceptance check.}
    \label{fig:posthoc_hist}
\end{figure}

Based on the observations outlined above, the warm-start analysis is re-evaluated using a post hoc acceptance check. 
Specifically, if an ART-generated trajectory satisfies the criteria $ -1.5 \leq {\mathcal{E}}_R \times 100 \leq 3.5$ and $ 0 \leq {\mathcal{E}}_C \leq 2$, it is accepted as a high-quality warm-start. 
If these conditions are not met, the CVX solution is employed as a fallback warm-start, resulting in a solution where $J_{\text{SCP,CVX}}-J_{\text{SCP,ART}} = 0$.
Figure \ref{fig:posthoc_hist} displays the histogram of the optimality gap when incorporating the post hoc acceptance check. 
This allows for a direct comparison with the bottom right subfigure of Figure \ref{fig:ws-analysis-chance-rtn}, highlighting the impact of the acceptance criteria on the performance of the warm-start strategies.

The impact of the post hoc acceptance check is observed in the results.
The minimum error of $J_{\text{SCP,CVX}}-J_{\text{SCP,ART}}$ approaches much closer to zero, effectively filtering out most warm-starting trajectories that would have negatively affected the SCP solution. 
However, the filtering process inadvertently eliminates some warm-starts that might have greatly improved the SCP's performance.
Nevertheless, the overall population is in the range of $-10 < J_{\text{SCP,CVX}}-J_{\text{SCP,ART}} < 500$ mm/s. 
Ultimately, the primary objective of the post hoc acceptance check---replacing low-quality ART-generated trajectories with fallback solutions---has been successfully achieved via simple, heuristic-based threshold definitions.
These results provide compelling evidence that runtime monitors can be developed to assess and enhance the performance of AI-based warm-starting techniques. 
Such monitors could focus on maximizing the retention of high-quality trajectories that lead to improved optimality while effectively filtering out low-quality solutions. 
Similar approaches represent an interesting and promising direction for future work.

% \subsection{Run-time Analysis on Flight GPU [TBD]}

%%%%%%%%%%%%%%%%%%%%%%%%%%%%%%%%%%%%%%%%%%%%%%%%%%%%%
\section{Conclusion}
%%%%%%%%%%%%%%%%%%%%%%%%%%%%%%%%%%%%%%%%%%%%%%%%%%%%%

This paper extends the capability of the Autonomous Rendezvous Transformer (ART), a Transformer-based warm-starting strategy, to address robust chance-constrained optimal control problems.
% to a non-convex chance-constrained optimal control problem. 
The proposed approach is evaluated in the context of robust fault-tolerant trajectory optimization for spacecraft rendezvous in Low Earth Orbit.
The experiments demonstrate that ART can effectively provide a high-quality robust trajectory as an initial guess for Sequential Convex Programming (SCP), resulting in faster convergence, improved optimality, and lower infeasibility rate.  
Furthermore, a post hoc acceptance criterion is developed to filter out low-quality generated trajectories, which is an essential prerequisite for deploying learning-based control systems in safety-critical environments. 
% It is demonstrated that simple criteria with a fallback solution can effectively prevent the convergence to degraded local optima. 
% This result motivates the research of the simple yet effective schemes of out-of-distribution detection (i.e., run-time monitors) and the development of alternative fallback initial guesses. 
There are several promising avenues for future work.
First, reducing the computational burden associated with transformer-based trajectory generation remains a priority, potentially leveraging model distillation and acceleration techniques. 
The run-time analysis of the flight-level GPU (e.g., NVIDIA Jetson series) should also be conducted to validate the real-time computation burden of the proposed method in a more realistic scenario.
Second, future research could explore incorporating failure scenarios beyond free-drift safety and expanding the applications to the cis-lunar domain, such as near-rectilinear Halo orbits.
Lastly, extending these methodologies to closed-loop control, for example through the integration with model predictive control techniques, is a compelling direction for advancing realistic autonomous operations in space. 

%%%%%%%%%%%%%%%%%%%%%%%%%%%%%%%%%%%%%%%%%%%%%%%%%%%%%%%%%%%%%%%%%%%%%%%%%%%%%%%%%%%%%%%%%%%%%%%%%

% \appendices{}              % note there is no {} to put a title. Each appendix has its own title
%%%%%%%%%%%%%%%%%%%%%%%%%%%%%%%%%%%%%%%%%%%%%%%%%%%%%%%%%%%%%%%%%%%%%%%%%%%%%%%%%%%%%%%%%%%%%%%%%
% For a single appendix, use the \appendix{} keyword and do not use the \section command.

\appendix{}\label{appendix:hyper} %: ART Hyperparameters
The Autonomous Rendezvous Transformer presented in this work is implemented in PyTorch \cite{paszke2017automatic} and builds off Huggingface's \code{transformers} library \cite{HuggFaceTransf}.
Specifically, Table \ref{tab_appendix:hyper} presents an overview of the hyperparameter settings used in this work.

\begin{table}[ht]
\centering
\caption{Hyperameters of ART architecture for the autonomous rendezvous experiments.}
\begingroup
\renewcommand*{\arraystretch}{1.25}
\begin{tabular}{l l}
    \hline 
    \hline
     Hyperparameter & Value \\
    \hline
     Number of layers & 6\\
     Number of attention heads & 6 \\
    Embedding dimension & 384 \\
     Batch size& 4 \\
    Context length $K$ & 100 \\
    Non-linearity & ReLU\\
    Dropout & 0.1\\
    Learning rate & 3e-5\\
    Grad norm clip & 1.0 \\
    Learning rate decay & None \\
    Gradient accumulation iters & 8 \\
    \hline
    \hline
    \end{tabular}%
    \label{tab_appendix:hyper}
    \endgroup
\end{table}

%%%%%%%%%%%%%%%%%%%%%%%%%%%%%%%%%%%%%%%%%%%%%%%%%%%%%%%%%%%%%%%%%%%%%%%%%%%%%%%%%%%%%%%%%%%%%%%%%%%%%%
\acknowledgments
This work is supported by Blue Origin (SPO \#299266) as an
Associate Member and Co-Founder of the Stanford’s Center of AEroSpace Autonomy Research (CAESAR). This article solely reflects the opinions and conclusions of its authors and not any Blue Origin entity.
Yuji Takubo acknowledges financial support from the Ezoe Memorial Recruit Foundation.

%%%%%%%%%%%%%%%%%%%%%%%%%%%%%%%%%%%%%%%%%%%%%%%%%%%%%%%%%%%%%%%%%%%%%%%%%%%%%%%%%%%%%%%%%%%%%%%%%%%%%%
\bibliographystyle{IEEEtran}
\bibliography{main} 
% \begin{thebibliography}{1}
% \bibitem{ITAR}
% U.S. Munitions List, Sections 38 and 47(7) of the Arms Export Control Act (22 U.S.C 2778 and 2794(7).
% \bibitem{AeroConf}
% Aerospace Conference Web site: \underline{www.aeroconf.org}.
% \end{thebibliography}

\clearpage

\thebiography
%% This biostyle allows you to insert your photo size 1in X 1.25in

\begin{biographywithpic}
{Yuji Takubo}{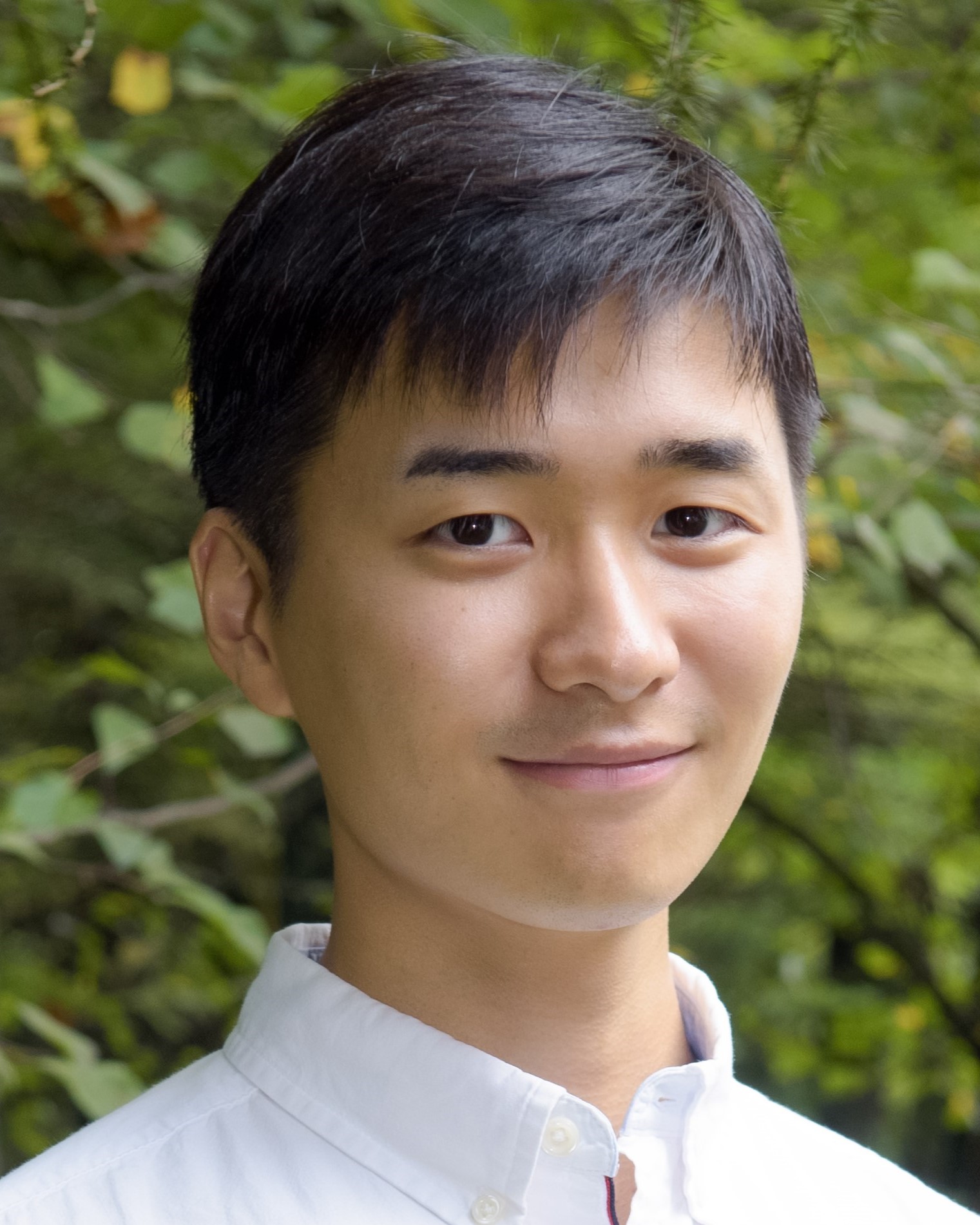}
Yuji Takubuo is a Ph.D. student in the Space Rendezvous Lab at Stanford University. He received the B.S. in Aerospace Engineering from the Georgia Institute of Technology (2023). 
He is a recipient of the NASA Jet Propulsion Laboratory Graduate Fellowship (2023) and the Ezoe Memorial Recruit Foundation Fellowship (2020-).  
His research focuses on optimization and astrodynamics, with topics ranging from multi-spacecraft GNC and interplanetary mission design to space logistics network optimization.
\end{biographywithpic}

\begin{biographywithpic}
{Tommaso Guffanti}{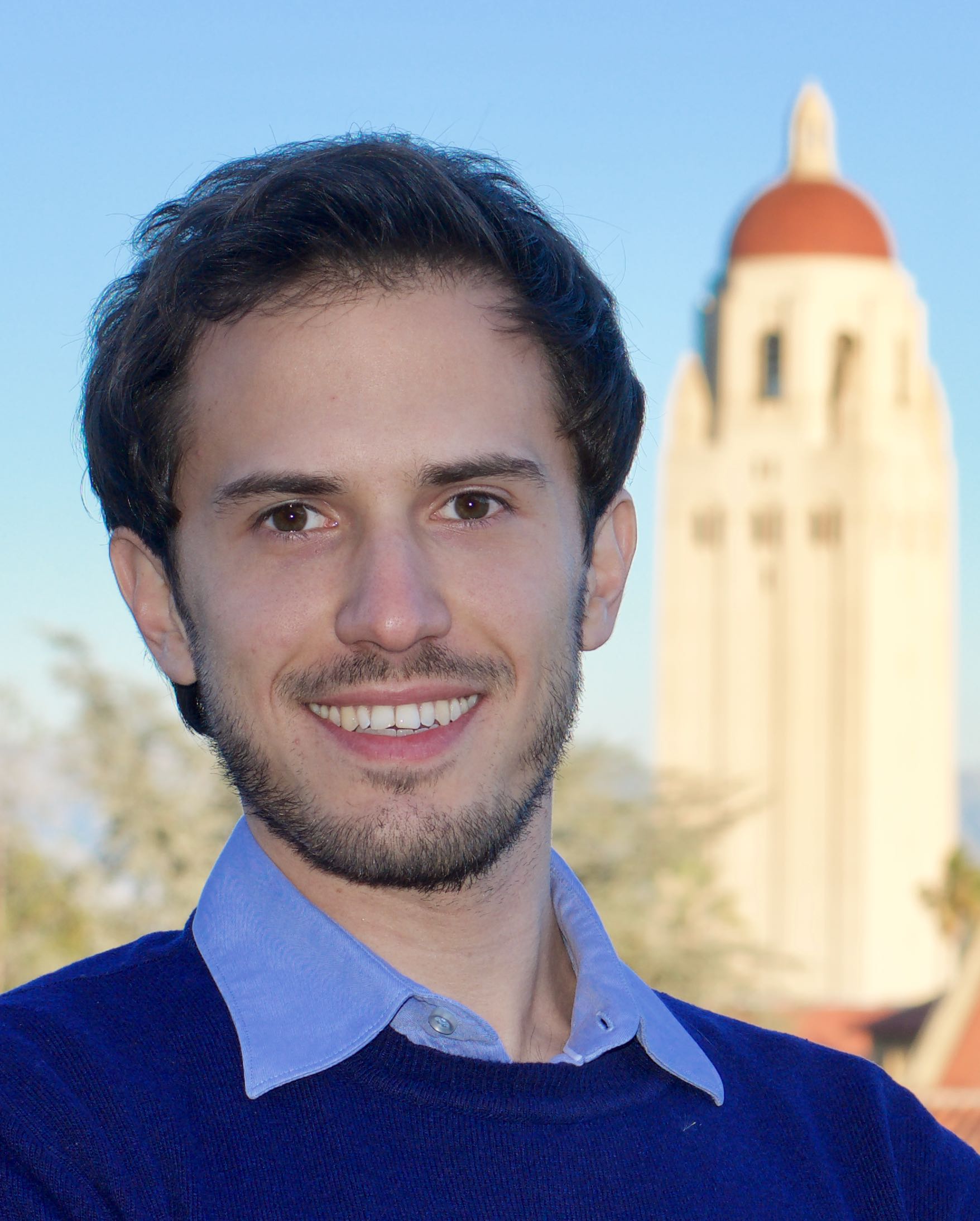}
Dr. Tommaso Guffanti is a Postdoctoral Scholar in the Space Rendezvous Lab at Stanford University. He received the B.S. and M.S. degrees in Aerospace engineering cum laude from Politecnico di Milano, and the Ph.D. degree in Aeronautics and Astronautics from Stanford University. Dr. Guffanti research contributions aim at developing cutting-edge guidance and control algorithms, and flight software to enable safe and autonomous functions and operations on-board space vehicles, in order to satisfy the requirements of the next generation of multi-spacecraft missions. During his doctorate and postdoctorate, he has been making research contributions in astrodynamics, safe and fault-tolerant multi-agent optimal control, and learning-based spacecraft motion planning and guidance for a variety of projects funded by national agencies and industry. He has over 15 scientific publications, including conference proceedings, and peer-reviewed journal articles. He has been awarded a doctoral Stanford Graduate Fellowship, a post-doctoral Stanford Center of Excellence for Aeronautics and Astronautics Scholarship, and a Stanford Emerging Technology Review Fellowship. He has been named excellent reviewer of the Journal of Guidance, Control, and Dynamics for three years.
\end{biographywithpic}

\begin{biographywithpic}
{Daniele Gammelli}{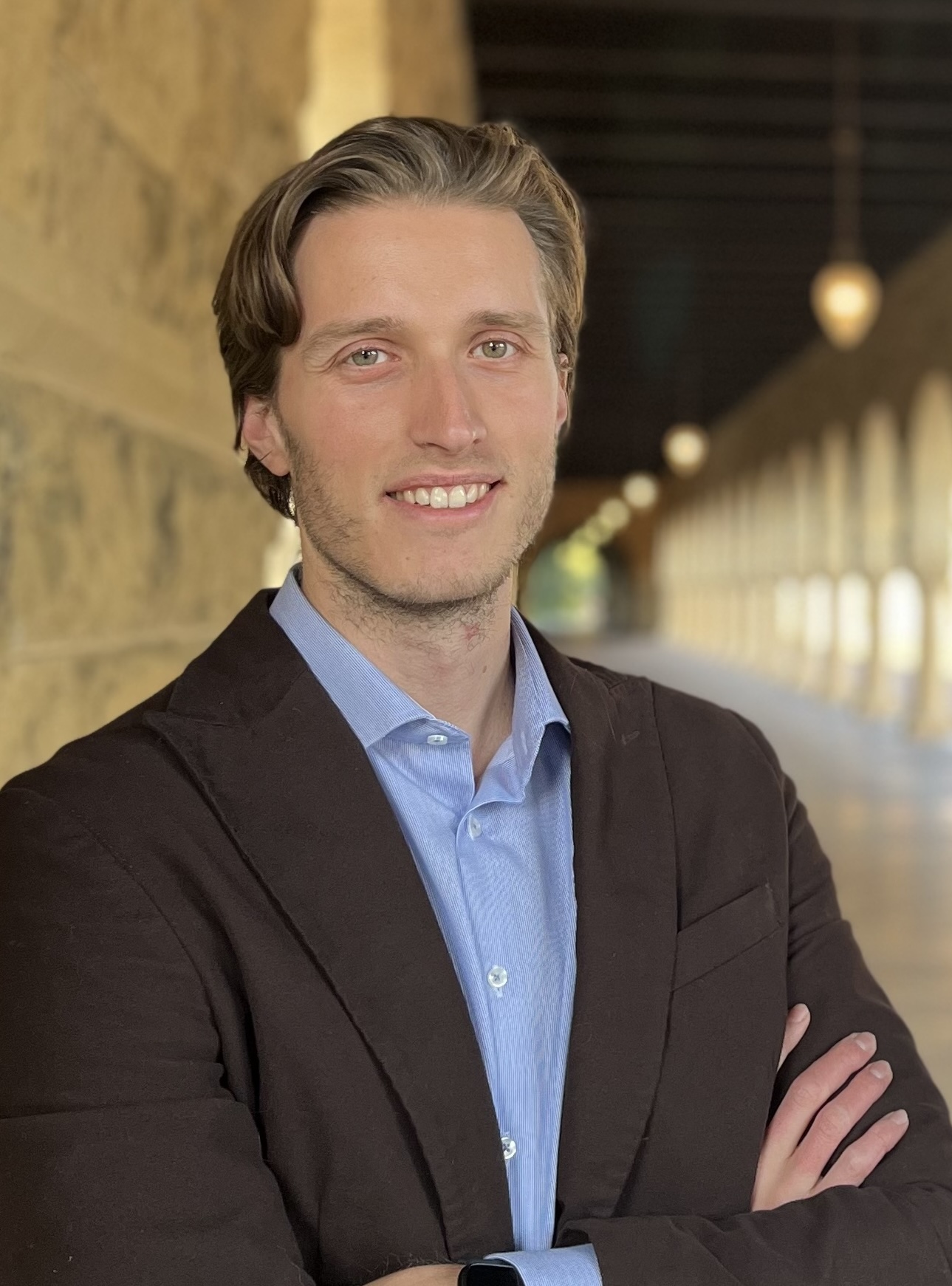}
Dr. Daniele Gammelli is a Postdoctoral Scholar in the Autonomous Systems Lab at Stanford University. He received the Ph.D. in Machine Learning and Mathematical Optimization at the Department of Technology, Management and Economics at the Technical University of Denmark. Dr. Gammelli’s research focuses on developing learning-based solutions that enable the deployment of future autonomous systems in complex environments, with an emphasis on large-scale robotic networks, aerospace systems, and future mobility systems. During his doctorate and postdoctorate career, Dr. Gammelli’s has been making research contributions in fundamental AI research, robotics, and its applications to network optimization and mobility systems. His research interests include deep reinforcement learning, generative models, graph neural networks, bayesian statistics, and control techniques leveraging these tools.
\end{biographywithpic}

\begin{biographywithpic}
{Marco Pavone}{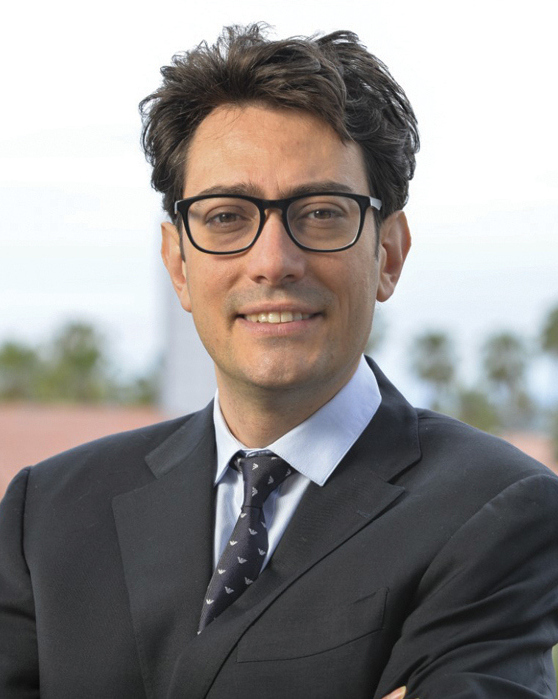}
Dr. Marco Pavone is an Associate Professor of Aeronautics and Astronautics at Stanford University, where he is the Director of the Autonomous Systems Laboratory and CoDirector of the Center for Automotive Research at Stanford. Before joining Stanford, he was a Research Technologist within the Robotics Section at the NASA Jet Propulsion Laboratory. He received a Ph.D. degree in Aeronautics and Astronautics from the Massachusetts Institute of Technology in 2010. His main research interests are in the development of methodologies for the analysis, design, and control of autonomous systems, with an emphasis on self-driving cars, autonomous aerospace vehicles, and future mobility systems. He is a recipient of a number of awards, including a Presidential Early Career Award for Scientists and Engineers from President Barack Obama, an Office of Naval Research Young Investigator Award, a National Science Foundation Early Career (CAREER) Award, a NASA Early Career Faculty Award, and an Early-Career Spotlight Award from the Robotics Science and Systems Foundation. He was identified by the American Society for Engineering Education (ASEE) as one of America’s 20 most highly promising investigators under the age of 40. His work has been recognized with best paper nominations or awards at the European Control Conference, at the IEEE International Conference on Intelligent Transportation Systems, at the Field and Service Robotics Conference, at the Robotics: Science and Systems Conference, at the ROBOCOMM Conference, and at NASA symposia. He is currently serving as an Associate Editor for the IEEE Control Systems Magazine. He is serving or has served on the advisory board of a number of autonomous driving start-ups (both small and multi-billion dollar ones), he routinely consults for major companies and financial institutions on the topic of autonomous systems, and is a venture partner for investments in AI-enabled robots.
\end{biographywithpic}

\begin{biographywithpic}
{Simone D'Amico}{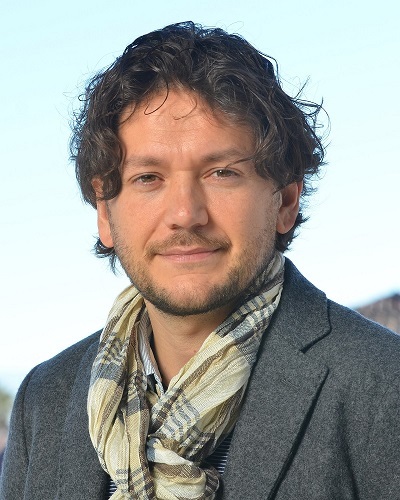}
Simone D’Amico is Associate Professor of Aeronautics and Astronautics (AA), W.M. Keck Faculty Scholar in the School of Engineering, and Professor of Geophysics (by Courtesy). He is the Founding Director of the Stanford Space Rendezvous Laboratory, Co-Director of the Center for AEroSpace Autonomy Research (CAESAR), and Director of the Undergraduate Program in Aerospace Engineering at Stanford. He has 20+ years of experience in research and development of autonomous spacecraft and distributed space systems. He developed the distributed Guidance, Navigation, and Control (GNC) system of several formation-flying and rendezvous missions and is currently the institutional PI of four autonomous satellite swarms funded by NASA (STARLING, STARI) and by NSF (VISORS, SWARM-EX) with one of them operational in orbit right now (Starling). Besides academia, Dr. D’Amico is in the Advisory Board of four space start-ups focusing on distributed space systems for future applications in SAR remote sensing, orbital lifetime prolongation, and space-based solar power. He was the recipient of several awards, most recently the 2024 NASA Ames Honor Award for the Starling mission, Best Paper Awards at IAF (2022), IEEE (2021), AIAA (2021), AAS (2019) conferences, and the M. Barry Carlton Award by IEEE (2020). He received the B.S. and M.S. degrees from Politecnico di Milano (2003) and the Ph.D. degree from Delft University of Technology (2010). 
\end{biographywithpic} 

\end{document}